\documentclass[12pt]{amsart}
\usepackage{amssymb,amsmath,amsthm}

\newtheorem{thm}{Theorem}
\newtheorem{lemma}[thm]{Lemma}

\newtheorem{prop}[thm]{Proposition}
\numberwithin{equation}{section} \theoremstyle{definition}

\newtheorem*{remark}{Remark}

\newcommand{\GL}{\mathrm{GL}}
\newcommand{\SL}{\mathrm{SL}}
\newcommand{\PSL}{\mathrm{PSL}}
\newcommand{\PGL}{\mathrm{PGL}}

\newcommand{\ve}{\varepsilon}
\newcommand{\Cay}{\mathrm{Cay}}
\newcommand{\Ind}{\mathrm{Ind}}
\newcommand{\mbfq}{\mathbb{F}_q}
\newcommand{\mbfQ}{\mathbb{E}}
\newcommand{\bW}{\overline{W}}
\newcommand{\ck}{\mathcal{K}}
\newcommand{\cl}{\mathcal{L}}
\newcommand{\cw}{\mathcal{W}}
\newcommand{\hf}{\hat{\mathbb{F}}}
\newcommand{\gd}{\delta}
\newcommand{\gG}{\Gamma}
\newcommand{\gh}{\theta}
\newcommand{\ghb}{\bar{\theta}}
\newcommand{\gl}{\lambda}
\newcommand{\gL}{\Lambda}
\newcommand{\go}{\omega}
\newcommand{\nf}{\infty}
\newcommand{\nr}{\text{N}}
\newcommand{\sd}{\sqrt{\delta}}
\newcommand{\tr}{\text{Tr}}
\newcommand{\tx}{^\times}

\begin{document}

\title[Ramanujan graphs on cosets of $\PGL_2(\mathbb F_q)$]{Ramanujan graphs on cosets of $\PGL_2(\mathbb F_q)$}
\author{Wen-Ching Winnie Li and Yotsanan Meemark}
\address{Wen-Ching Winnie Li\\ Department of Mathematics\\ The Pennsylvania State University\\ University Park, PA 16802}
\email{\tt wli@math.psu.edu}
\address{Yotsanan Meemark\\ Department of Mathematics\\ The Pennsylvania State University\\ University Park, PA 16802}
\email{\tt meemark@math.psu.edu}
\thanks{The research of the first author is supported in part by an NSA
grant MDA904-03-1-0069.}


\begin{abstract}

In this paper we study Cayley graphs on  $\PGL_2(\mathbb F_q)$ mod
the unipotent subgroup, the split and nonsplit tori, respectively.
Using the Kirillov models of the representations of
$\PGL_2(\mathbb F_q)$ of degree greater than one, we obtain
explicit eigenvalues of these graphs and the corresponding
eigenfunctions. Character sum estimates are then used to conclude
that two types of the graphs are Ramanujan, while the third is
almost Ramanujan. The graphs arising from the nonsplit torus were
previously studied by Terras et al. We give a different approach
here.
\end{abstract}
\maketitle

\section{Introduction}
A finite $k$-regular graph is called Ramanujan if its eigenvalues
other than $\pm k$, called nontrivial eigenvalues, have absolute
values at most $2\sqrt {k-1}$. Such graphs are good expanders and
have broad applications in computer science. The first systematic
explicit construction of an infinite family of $k$-regular
Ramanujan graphs is given independently by Margulis \cite{Ma88}
and Lubotzky-Phillips-Sarnak \cite{LPS88} for $k = p+1$ with $p$ a
prime; their graphs are based on quaternion groups over $\mathbb
Q$, and the nontrivial eigenvalues of these graphs can be
interpreted as the eigenvalues of the Hecke operators on classical
cusp forms of weight 2 so that the eigenvalue bound follows from
the deep property that the Ramanujan conjecture holds for these
cusp forms, established by Eichler \cite{Ei57} and Shimura
\cite{Sh59}. The Ramanujan conjecture for cusp forms of $\GL_2$
over a function field is proved by Drinfeld \cite{Dr88}, hence the
same method gives rise to infinite families of $(q + 1)$-regular
Ramanujan graphs, where $q$ is a prime power. This is done by
Morgenstern in \cite{Mo94}.

On the other hand, there are explicit constructions of
$(q+1)$-regular Ramanujan graphs for $q$ a prime power whose
nontrivial eigenvalues are expressed as character sums, which are
shown to be bounded by $2\sqrt q$, as a consequence of the Riemann
hypothesis for curves over finite fields. Such examples include
Terras graphs \cite{CP93} based on cosets of $\PGL_2(\mbfq)$ and
norm graphs in \cite{L92} based on a finite field of $q^2$
elements. Later it is shown in \cite{L99} that these two kinds of
graphs are in fact quotient graphs of Morgenstern graphs. This
connection leads to very interesting relations between character
sums and cusp forms for $\GL_2$ over function fields, studied in
detail in \cite{CL03} and \cite{CL04}. In particular, one obtains
cusp forms whose Fourier coefficients are given by eigenvalues of
the Terras/norm graphs in a systematic manner.

In this paper we revisit Terras graphs and investigate two other
types of graphs based on cosets of $G = \PGL_2(\mbfq)$ using a
uniform method explained below. Up to conjugation, $G$ contains
three types of abelian subgroups: the unipotent subgroup
\begin{eqnarray*}
U = \{ \left(\begin{smallmatrix}
  1 & x \\
  0 & 1 \\
\end{smallmatrix}\right) : x \in \mbfq \},
\end{eqnarray*}
which is a group of order $q$, the split torus
\begin{eqnarray*}
A = \{ \left(\begin{smallmatrix}
  y & 0 \\
  0 & 1 \\
\end{smallmatrix}\right): y \in \mathbb F_q^\times \},
\end{eqnarray*}
which is cyclic of order $q-1$, and the nonsplit torus $K$, which
is an embedded image of $\mathbb F_{q^2}^\times/\mathbb
F_q^\times$ in $G$ as a cyclic subgroup of order $q+1$. It is
well-known that $G = UAK$. Denote by $H$ one of these three
subgroups. For a double coset $HsH$ which is its own inverse
(i.e., symmetric) and which is the disjoint union of $|H|$ right
$H$-cosets, consider the Cayley graph $X_{HsH} = \Cay(G/H,
HsH/H)$, called an $H$-graph. It is undirected and $|H|$-regular.
(When $H = K$, this is a Terras graph.) We shall prove
\medskip

\noindent {\bf Main Theorem} (a) (cf. Terras \cite{Ta99}, p.357) {\it The nontrivial eigenvalues of $X_{KsK}$ have 
absolute values at most $2\sqrt q$. Hence the $K$-graphs are $(q+1)$-regular
Ramanujan graphs.} 

(b) {\it The nontrivial eigenvalues of $X_{UsU}$ are $\pm 1$ and $\pm \sqrt q$.
Thus the $U$-graphs are $q$-regular Ramanujan graphs.}

(c) {\it The nontrivial eigenvalues of $X_{AsA}$ have 
absolute values at most $2\sqrt q$. Thus the $A$-graphs, being 
$(q-1)$-regular, are almost Ramanujan.}
\medskip


Like representations of $p$-adic groups, the irreducible
representations of $G$ of degree greater than one also have a
Kirillov model, in which the actions of $U$ and $A$ are standard,
and the representations are distinguished by the action of the
Weyl element $w = \left(\begin{smallmatrix}
  0 & 1 \\
  -1 & 0 \\
\end{smallmatrix}\right)$. This is studied in \cite{LSA83}. Our
approach is to use a Kirillov model to explicitly find functions
in the space of each irreducible representation of $G$ which are
right $H$-invariant. This space has dimension at most 2. We then
determine the eigenvalues and eigenfunctions of the adjacency
operator of the graph $X_{HsH}$. For Terras graphs, the
eigenvalues are obtained in \cite{AC92}, \cite{CP93} by computing
the traces of the irreducible representations, expressed as
character sums over $\mbfq$, and then estimated. In our approach
we obtain the same expression for eigenvalues arising from
nondiscrete series representations, but an eigenvalue arising from
a discrete series representation is expressed differently, namely,
as the average of $q+1$ character sums over $\mathbb F_{q^2}$, one
for each $K$-coset. Each character sum is associated with some
id\`{e}le class character of the rational function field $\mathbb
F_q(T)$ using the results obtained in Chapter 6 of \cite{L96} and
\cite{L05}, and then shown to have absolute value at most $2\sqrt
q$ as a consequence of the Riemann hypothesis for curves.

It is worth pointing out that the $U$- graph $X_{UwU}$ has two
connected components, one of which can be identified with the
Cayley graph $\Cay(\PSL_2(\mbfq)/U, UwU/U)$. A suitable quotient
of this graph may be interpreted as a graph on the cusps of a
certain principal congruence subgroup of the Drinfeld modular
group $\GL_2(\mathbb F_p[T])$. This kind of graph is first
obtained by Gunnells \cite{G04} for principal congruence subgroups
$\Gamma(p)$ of $\SL_2(\mathbb Z)$. A generalization of this graph
from $p$ to prime power $q$ with $q \equiv 1 \pmod 4$ is given in
\cite{De04}. Both approaches rely on analyzing the graph
structure, while ours is purely representation-theoretic.

The paper is organized as follows. The representation theory,
including the Kirillov models, is reviewed in section 2. Sections
3, 4, and 5 are devoted to the $K$-, $U$-, and $A$-graphs,
respectively. In each case, using Kirillov models, we determine
the eigenvalues and compute the corresponding eigenfunctions; then
character sum estimates are employed to bound the eigenvalues.

For convenience, the characteristic of $\mathbb F_q$ is assumed to
be odd throughout the paper. Similar results are expected to hold
for even characteristic.

\bigskip

\section{Representations of $\PGL_2(\mathbb{F}_q)$}

For brevity, write $\mathbb{F}$ for the finite field with $q$
elements and $\mathbb{E}$ for its quadratic extension. The
unipotent subgroup $U$ acts on the space $\cl(G)$ of
complex-valued functions  via left translations so that the space
decomposes as
\[
\cl(G) = \bigoplus_{\psi \in \hat{\mathbb{F}}} \cl_\psi (G)
\]
where $\cl_\psi (G) = \{ f: G \to \mathbb{C} : f
\left(\left(\begin{smallmatrix}
  1 & x \\
  0 & 1 \\
\end{smallmatrix}\right)g\right) = \psi(x)f(g)$ for all $g \in
G\}$.

The irreducible representations of $G$ are identified with the
irreducible representations of $\GL_2(\mathbb{F})$ with trivial
central character. Such representations are studied in the
literature in detail (cf. \cite{PS83}); those with degree greater
than one are classified into three categories: principal series,
Steinberg, and discrete series representations. Fix a nontrivial
additive character $\psi$ of $\mathbb{F}$. Each irreducible
representation $\pi$ of $G$ of degree greater than 1 has a
Kirillov model $\ck_\psi (\pi)$ and a Whittaker model
$\cw_\psi(\pi)$. According to $\pi$ being a discrete series,
Steinberg, or principal series representation, the space $\ck_\psi
(\pi)$ is spanned by $\hat{\mathbb{F}}^\times$ (the multiplicative
characters of $\mathbb{F}$), $\hat{\mathbb{F}}^\times \cup \{ D_0
\}$, or $\hat{\mathbb{F}}^\times \cup \{ D_0 \} \cup \{
D_\infty\}$ respectively, where $D_0$ (resp. $D_\infty$) denotes
the Dirac function at $0$ (resp. $\infty$). The action of $U$ on
$\hf^\times$ is given in terms of $\psi$, and the action of $UA$
on $\hf^\times$ is the same for all representations. The
representation $\pi$ is characterized by the action of $w =
\left(\begin{smallmatrix}
  0 & 1 \\
  -1 & 0 \\
\end{smallmatrix}\right)$. In \cite{LSA83}, to each $\pi$ of degree greater
than 1 a family of Gauss sums $\ve(\pi, \chi, \psi)$, where $\chi
\in \hf^\times$, is attached, with which the action of $\pi(w)$ is
described. The details are as follows.

A discrete series representation $\pi = \pi_\Lambda$ arises from a
character $\Lambda$ of $\mathbb{E}^\times$ such that $\Lambda$ is
trivial on $\mathbb{F}^\times$ and $\Lambda \ne \mu \circ \nr$ for
all $\mu \in \hf^\times$. Here $\nr$ denotes the norm map from
$\mathbb{E}$ to $\mathbb{F}$, The last condition on $\Lambda$ is
equivalent to $\Lambda^{q+1} = 1, \Lambda^2 \ne 1.$ Consequently,
the inverse of $\Lambda$ is $\bar{\Lambda} = \Lambda^q$. We
associate $\pi$ with
\[
\ve(\pi, \chi, \psi) = - \Gamma(\Lambda\chi \circ \text{N}, \psi
\circ \tr) = - \sum_{z \in \mathbb{E}^\times} \Lambda(z) \chi(\nr
\, z)\psi(\tr \, z), \quad \chi \in \hf^\times.
\]
Here $\tr$ is the trace map from $\mathbb{E}$ to $\mathbb{F}$.
Observe that $\ve(\pi_\Lambda, \chi, \psi) =
\ve(\pi_{\bar{\Lambda}}, \chi, \psi)$.

Characters $\mu$ of $\hf^\times$ give rise to principal series
representations $\pi_\mu$ if $\mu^2 \ne 1$, and Steinberg
representations $\pi_\mu$ if $\mu^2 = 1$. These are all
nondiscrete series representations of $G$ of degree greater than
1. For $\pi = \pi_\mu$ arising from a character $\mu$ of
$\hf^\times$, we associate
\[
\ve(\pi, \chi, \psi) = \gG(\mu\chi, \psi) \gG(\mu^{-1}\chi, \psi),
\quad \chi \in \hf^\times,
\]
where $\gG(\xi, \psi) = \sum_{x \in \mathbb{F}^\times}
\xi(x)\psi(x)$. Observe that $\ve(\pi_{\mu}, \chi, \psi) =
\ve(\pi_{\mu^{-1}}, \chi, \psi)$.

Using $\ve(\pi, \chi, \psi)$, we can describe the representation
$\pi$ on its Kirillov model $\ck_\psi(\pi)$ by giving the action
of the generators
\[
h_r = \left(\begin{smallmatrix}
  r & 0 \\
  0 & 1 \\
\end{smallmatrix}\right)\;(r \in \mathbb{F}^\times), u_s = \left(\begin{smallmatrix}
  1 & s \\
  0 & 1 \\
\end{smallmatrix}\right)\;(s \in \mathbb{F}^\times), w = \left(\begin{smallmatrix}
  0 & 1 \\
  -1 & 0 \\
\end{smallmatrix}\right)
\]
of $G$ as follows:
\[
\begin{aligned}
\pi(h_r)\alpha &= \alpha(r)\alpha, &&\alpha \in \hf^\times, \\
\pi(h_r)D_0 &= \mu(r)D_0, \pi(h_r)D_\infty = \mu^{-1}(r)D_\infty
&&\text{(when applicable),}\\
\pi(u_s)\alpha &= (q-1)^{-1} \sum_{\beta \in \hf^\times}
\beta\alpha^{-1}(s)\gG(\alpha\beta^{-1}, \psi)\beta, &&\alpha \in
\hf^\times,
\\
\pi(u_s)D_0 &= D_0, \pi(u_s)D_\infty = D_\infty &&\text{(when applicable),}\\
\pi(w)\alpha &= q^{-1}\ve(\pi, \alpha, \psi)\alpha^{-1} +
e(\alpha), &&\alpha
\in \hf^\times, \\
\pi(w)D_0 &= -q^{-1}\ve(\pi, \mu^{-1}, \psi)(\mu^{-1}+D_0)
&&\text{if $\pi$
is a Steinberg representation}, \\
&=-q^{-1}\ve(\pi, \mu, \psi)(\mu^{-1}+D_\infty) &&\text{if $\pi$
is a principal series}, \\
\pi(w)D_\infty &= -q^{-1}\ve(\pi, \mu^{-1}, \psi)(\mu+D_0)
&&\text{if $\pi$
is a principal series}, \\
\end{aligned}
\]
where
\[
\begin{aligned}
e(\alpha) &=0 &&\text{if $\pi$ is a discrete series},\\
&= - q^{-1}\ve(\pi, \alpha, \psi)(q^2 -1)\delta_{\alpha\mu, 1}D_0
&&\text{if
$\pi$ is a Steinberg representation}, \\
&= - q^{-1}\ve(\pi, \alpha, \psi)(q-1)(\delta_{\alpha\mu,1}D_0
+\delta_{\alpha\mu^{-1}, 1}D_\infty) &&\text{if $\pi$ is a
principal series}.
\end{aligned}
\]
Here $\gd_{\chi, 1}$ is the Kronecker symbol, which is equal to
$1$ if $\chi = 1$, the trivial character, and $0$ otherwise. It
was shown in \cite{LSA83} that the relations on these generators
are preserved, resulting from the identities satisfied by the
Gauss sums with the main one called the Barnes' identity. The
representations are characterized by the attached $\ve$-factors.
In particular, $\pi_{\gL}$ and $\pi_{\bar{\gL}}$ are equivalent,
and so are $\pi_\mu$ and $\pi_{\mu^{-1}}$.

In conclusion, there are $\frac{q-1}{2}$ (nonequivalent) discrete
series representations $\pi_\gL$, each of degree $q-1$; there are
two Steinberg representations $\pi_\mu$ from $\mu$ with $\mu^2 =
1$, each of degree $q$; and there are $\frac{q-3}{2}$ principal
series representations $\pi_\mu$ with $\mu^2 \ne 1$, each of
degree $q+1$. Together with the two degree one  representations
given by $\mu \circ \det$, $\mu^2 = 1$, this gives the complete
list of the irreducible representations of $G$. Each irreducible
representation occurs in $\cl(G)$ with multiplicity equal to its
degree.

The Whittaker model $\cw_\psi(\pi)$ consists of functions on $G$
obtained from the Kirillov model $\ck_\psi(\pi)$ in the following
way: for $v \in \ck_\psi(\pi)$, define $W_v \in \cw_\psi(\pi)$ via
\[
W_v(g) = (\pi(g)v)(1), \quad g \in G.
\]
The representation $\pi$ of $G$ on $\cw_\psi(\pi)$ is by right
translations. Because of the action of $U$ on $\ck_\psi(\pi)$, all
functions are contained in $\cl_\psi(G)$.

We proceed to discuss how each space $\cl_\psi(G)$ decomposes. Fix
a nontrivial additive character $\psi$ of $\mathbb{F}$. We can
describe all characters of $\mathbb{F}$ as $\psi^a$, $a \in
\mathbb{F}$, where $\psi^a(x) = \psi(ax)$ for $x$ in $\mathbb{F}$.
Observe that for each function $f \in \cl_\psi(G)$, the new
function $f_a(g) := f(\left(\begin{smallmatrix}
  a & 0 \\
  0 & 1 \\
\end{smallmatrix}\right)g)$ lies in $\cl_{\psi^a}(G)$ and $f \mapsto
f_a$ is an isomorphism between the two $G$-modules $\cl_\psi(G)$
and $\cl_{\psi^a}(G)$. As noted above, the Whittaker model
$\cw_\psi(\pi)$ of each irreducible $\pi$ of degree greater than 1
is contained in $\cl_\psi(G)$ for $\psi$ nontrivial, by counting
dimension we find
\[
\cl_{\psi^a}(G) = \bigoplus_{\pi, \deg \pi >1} \cw_{\psi^a}(\pi),
\quad \text{for $a \ne 0$}.
\]
We then conclude from checking the multiplicities that
$\cl_{\psi^0} (G)$ contains no discrete series representations;
 each principal series representation occurs there twice,
and each Steinberg representation and each degree 1 representation
occur once. We record this in

\begin{prop} \label{rep}
\begin{itemize}
\item[(1)] For nontrivial $\psi$'s, the spaces $\cl_\psi(G)$ are
isomorphic and each irreducible representation of $G$ of degree
greater than $1$ occurs in $\cl_\psi(G)$ exactly once.

\item[(2)] For $\psi = \psi^0$ trivial, in $\cl_{\psi^0}(G)$ each
principal series representation occurs with multiplicity two, each
Steinberg representation occurs with multiplicity one, as does
each one dimensional representation.
\end{itemize}
\end{prop}

\begin{remark}
For each character $\mu$ of $\mathbb{F}^\times$, denote by
\[
\Ind\mu = \{f:G \to \mathbb{C}: f \left(\left(\begin{smallmatrix}
  a & x \\
  0 & 1 \\
\end{smallmatrix}\right)g\right) = \mu(a)f(g) \; \text{for all} \; g \in
G\}.
\]
The group $G$ acts on $\Ind \mu$ by right translations. It is
well-known that when $\mu^2 \ne 1$, this representation is the
principal series representation $\pi_\mu$, and it is isomorphic to
$\Ind \mu^{-1}$. When $\mu^2 = 1$, this representation has two
irreducible constituents, $\mu \circ \det$ and the Steinberg
representation $\pi_\mu$. Obviously,
\[
\cl_{\psi^0} (G) = \bigoplus_{\mu \in \hf^\times} \Ind \mu.
\]
\end{remark}

Let $H$ be a subgroup of $G$. Then the space of functions on $G$
right invariant by $H$ has a similar decomposition
\[
\cl(G/H) = \bigoplus_\psi \cl_\psi(G/H) = \bigoplus_\psi
\bigoplus_\pi m_{\psi}(\pi) \cl_\psi(\pi, G/H),
\]
where $\pi$ runs through irreducible representations of $G$,
$m_\psi(\pi)$ is the multiplicity of $\pi$ in $\cl_\psi(G)$ as
described in the proposition above, and $\cl_\psi (\pi, G/H)$
consists of the right $H$-invariant functions in the space of
$\pi$ in $\cl_\psi(G)$.

Let $s$ be an element in $G$. Write $HsH = x_1H \cup \dots \cup
x_kH$ as a disjoint union of $k$ right $H$ cosets. Define an
operator $T_{HsH}$ on $\cl(G/H)$ by sending $f \in \cl(G/H)$ to
\[
(T_{HsH}f)(xH) = \sum_{i=1}^k f(xx_iH).
\]
Clearly $T_{HsH}$ preserves each space $\cl_\psi(\pi, G/H)$. When
$HsH = Hs^{-1}H$, we define an undirected Cayley graph $X_{HsH} =
\Cay(G/H, HsH/H)$, called an $H$-graph, whose adjacency matrix may
be identified with the operator $T_{HsH}$. We shall take $H=K, U$
and $A$, respectively, and study the eigenfunctions and the
eigenvalues of $T_{HsH}$.

\bigskip

\section{The $K$-graphs}

In this section we take $H=K$. Since $G=UAK$, for each additive
character $\psi$ of $\mathbb{F}$, the space $\cl_\psi (G/K)$ is
$(q-1)$-dimensional. As remarked before, when $\psi = \psi^0$ is
the trivial character,
\[
\cl_{\psi^0}(G) = \bigoplus_{\mu \in \hf^\times} \Ind \mu.
\]
One sees immediately from the definition of $\Ind \mu$ that the
right $K$-invariant space of $\Ind \mu$ is 1-dimensional,
generated by
\[
f_\mu\left(\left(\begin{smallmatrix}
  y & x \\
  0 & 1 \\
\end{smallmatrix}\right)k\right) = \mu(y) \quad \text{for all $y \in
\mathbb{F}^\times, x \in \mathbb{F}$ and $k \in K$},
\]
and the $f_\mu$'s form a basis of $\cl_{\psi^0}(G)$. This shows
that for a principal series representation $\pi_\mu$, the space of
right $K$-invariant vectors in any of its model is 1-dimensional
for any $\psi$. When $\mu$ is the quadratic character, the
1-dimensional representation $\mu \circ \det$ does not contain
nontrivial right $K$-invariant vectors since not all elements in
$K$ have square determinant. Therefore $f_\mu$ belongs to the
Steinberg representation $\pi_\mu$, and we arrive at the same
conclusion that $\cl_\psi(\pi_\mu, G/K)$ is 1-dimensional. When
$\mu$ is the trivial character, the function $f_\mu$ lies in the
space of the 1-dimensional representation $\mu \circ \det$, and
hence $\cl_\psi(\pi_\mu, G/K)$ is 0-dimensional for the Steinberg
representation $\pi_\mu$ for $\mu$ trivial. We record this in

\begin{prop}
Let $\mu$ be a character of $\mathbb{F}^\times$. For any additive
character $\psi$ of $\mathbb{F}$, the space $\cl_\psi (\pi_\mu,
G/K)$ is $1$-dimensional if $\mu \ne 1$, and $0$-dimensional if
$\mu = 1$.
\end{prop}

Our next goal is to show that, for $\psi$ nontrivial,
$\cl_\psi(\pi_\gL, G/K)$ is 1-dimensional for each discrete series
representation $\pi_\gL$. Since $\cl_\psi(G/K)$ is
$(q-1)$-dimensional and it contains a $\frac{q-1}{2}$-dimensional
subspace $\bigoplus_{\mu \ne 1} \cl_\psi(\pi_\mu, G/K)$, it
suffices to show $\dim \cl_\psi (\pi_\gL, G/K) \ge 1$ as there are
$\frac{q-1}{2}$ discrete series representations.

{\it Fix once and for all a nonsquare $\delta$ in $\mathbb{F}$} so
that $\mathbb{E} = \mathbb{F}(\sqrt{\delta})$. We imbed
$\mathbb{E}^\times$ in $\GL_2(\mathbb{F})$ as
$\{\left(\begin{smallmatrix}
  b & a\delta \\
  a & b \\
\end{smallmatrix}\right) : a,b \in \mathbb{F} \; \text{not both zero} \}.$
Consequently, the elements in $K$ are represented by
\[
\left(\begin{smallmatrix}
  b & \delta \\
  1 & b \\
\end{smallmatrix}\right) = \left(\begin{smallmatrix}
  1 & b \\
  0 & 1 \\
\end{smallmatrix}\right)\left(\begin{smallmatrix}
  0 & 1 \\
  -1 & 0 \\
\end{smallmatrix}\right)\left(\begin{smallmatrix}
  1 & b(b^2-\delta)^{-1} \\
  0 & 1 \\
\end{smallmatrix}\right)\left(\begin{smallmatrix}
  1 & 0 \\
  0 & b^2-\delta \\
\end{smallmatrix}\right)\left(\begin{smallmatrix}
  -1 & 0 \\
  0 & -1 \\
\end{smallmatrix}\right) \,\ \text {with} ~~b \in \mathbb{F},  \quad \text{and} \quad
\left(\begin{smallmatrix}
  1 & 0 \\
  0 & 1 \\
\end{smallmatrix}\right).
\]
Fix a nontrivial additive character $\psi$ of $\mathbb{F}$ and a
character $\gL$ of $\mathbb{E}^\times$ with $\gL^{q+1}=1$ and
$\gL^2 \ne 1$. Let $\pi = \pi_\gL$ be the associated discrete
series representation with the attached $\ve$-factor $\ve(\pi,
\chi, \psi)$. According to the description in \S 2 of the action
of $\pi$ on its Kirillov model $\ck_\psi(\pi)$, we obtain, for
$\theta \in \hf^\times$ and $b \in \mathbb{F}^\times$,
\[
{\small\begin{split} \pi\left(\left(\begin{smallmatrix}
  b & \delta \\
  1 & b \\
\end{smallmatrix}\right)\right)\theta =
(q-1&)^{-2}q^{-1}\bar{\theta}(b)\sum_{\beta \in \hf^\times}
\beta(b(b^2 - \delta)^{-1}) \gG(\theta\bar{\beta}, \psi)\ve(\pi,
\beta, \psi) \sum_{\gamma \in \hf^\times} (\gamma\beta)(b)
 \gG(\bar{\gamma}\bar{\beta}, \psi)
\gamma, \\
&\pi\left(\left(\begin{smallmatrix}
  0 & \delta \\
  1 & 0 \\
\end{smallmatrix}\right)\right)\theta = \bar{\theta}(-\delta) \ve(\pi, \theta,
\psi)q^{-1}\bar{\theta} \quad \text{and} \quad
\pi\left(\left(\begin{smallmatrix}
  1 & 0 \\
  0 & 1 \\
\end{smallmatrix}\right)\right)\theta = \theta.
\end{split}}
\]
Then $v_\theta := \sum_{b \in \mathbb{F}}
\pi\left(\left(\begin{smallmatrix}
  b & \delta \\
  1 & b \\
\end{smallmatrix}\right)\right)\theta + \theta$ is invariant by
$\pi(K)$. Put
\[
{\small\begin{aligned}
v_{\pi, \psi} &= \sum_{\gh \in \hf^\times} v_\gh \\
&= \sum_{\gh} \sum_{b \in \hf^\times} (q-1)^{-2}q^{-1} \ghb(b)
\sum_{\beta \in \hf^\times}
\beta(b(b^2-\delta)^{-1})\gG(\gh\bar{\beta}, \psi) \ve(\pi, \beta,
\psi) \sum_{\gamma \in \hf^\times}(\gamma\beta)(b)\gG(\bar{\gamma}
\bar{\beta}, \psi)\gamma  \\
& \quad + \sum_{\gh} \ghb(-\delta)\ve(\pi, \gh, \psi) q^{-1} \ghb
+ \sum_{\gh} \gh.
\end{aligned}}
\]
Using
\[
\sum_{\gh \in \mathbb{F}\tx} \ghb(b)\gG(\gh\beta, \psi) =
\sum_{\gh} \ghb(b) \sum_{x \in \mathbb{F}^\times}
(\gh\bar{\beta})(x)\psi(x) = (q-1)\bar{\beta}(b)\psi(b),
\] we
rewrite $v_{\pi, \psi}$ as
\[
{\small\begin{aligned} v_{\pi, \psi} &= \sum_{b \in
\mathbb{F}^\times} (q-1)^{-1}q^{-1} \sum_{\beta \in \hf^\times}
\bar{\beta}(b^2 - \delta) \psi(b) \ve(\pi, \beta, \psi)
\sum_{\gamma \in \hf^\times}
(\gamma\beta)(b)\gG(\bar{\gamma}\bar{\beta},
\psi)\gamma \\
& \quad + \sum_{\gh} \ghb(-\delta)\ve(\pi, \gh, \psi) q^{-1} \ghb
+ \sum_{\gh} \gh.
\end{aligned}}
\]
Recall that $\ve(\pi, \chi, \psi) = - \Gamma(\Lambda \chi\circ
\nr, \psi \circ \tr) = - \sum_{z \in \mathbb{E}^\times}
\Lambda(z)\chi(\nr ~z)\psi(\tr ~z)$. Summing over $\beta$ further
simplifies the expression of $v_{\pi, \psi}$ as
\[
{\small\begin{split}v_{\pi, \psi} = -q^{-1}\sum_{\gamma \in
\hf^\times} \sum_{b \in \mathbb{F}} \sum_{z \in \mathbb{E}^\times}
\gG(z) \bar{\gamma}(\nr\,z)\psi(\tr\,z)\gamma(b^2-\delta)
\psi((\nr z)b(b^2-\delta)^{-1})\psi(b)\gamma + \sum_\gamma \gamma.
\end{split}}
\]
Let $W_{\gL, \psi}' = W_{v_{\pi, \psi}}$, that is, for $g \in G$,
$W_{\gL, \psi}(g) = (\pi(g)v_{\pi,\psi})(1)$. Then $W_{\gL,
\psi}'$ lies in the Whittaker model $\cw_\psi(\pi)$ and it is
right $K$-invariant. Hence $W_{\gL, \psi}'$ is determined by its
values on the split torus $A$, which we now compute. For $y \in
\mathbb{F}^\times$,
\[
{\small\begin{aligned}
&W_{\gL,\psi}'\left(\left(\begin{smallmatrix}
  y & 0 \\
  0 & 1 \\
\end{smallmatrix}\right)\right) = (\pi\left(\left(\begin{smallmatrix}
  y & 0 \\
  0 & 1 \\
\end{smallmatrix}\right)\right)v_{\pi,\psi})(1) \\
&=-q^{-1}\sum_{\gamma \in \hf^\times} \sum_{b \in \mathbb{F}}
\sum_{z \in \mathbb{E}^\times} \gG(z)
\bar{\gamma}(\nr\,z)\psi(\tr\,z)\gamma(b^2-\delta) \psi((\nr
z)b(b^2-\delta)^{-1})\psi(b)\gamma(y) + \sum_\gamma \gamma(y) \\
&=-(q-1)q^{-1}\sum_{b \in \mathbb{F}} \sum_{z \in
\mathbb{E}^\times, \atop \nr\,z = y(b^2-\delta)} \gL(z)\psi(\tr\,
z) \psi((\nr \, z) b(b^2-\delta)^{-1})\psi(b) + (q-1)\delta_{y,1}.
\end{aligned}}
\]
Here $\gd_{y,1}$ is equal to $1$ if $y=1$ and $0$ otherwise.
Replacing the variable $z$ by $z(b+\delta)$, we rewrite the above
as
\[
{\small\begin{split} W_{\gL,\psi}'\left(\left(\begin{smallmatrix}
  y & 0 \\
  0 & 1 \\
\end{smallmatrix}\right)\right) = -(q-1)q^{-1}\sum_{b \in \mathbb{F}}
\sum_{z \in \mathbb{E}^\times, \atop \nr\,z =y}
\gL(z(b+\sqrt{\delta}))\psi(\tr(z(b+\sqrt{\delta})))\psi(b(y+1)) +
(q-1)\delta_{y,1}.
\end{split}}
\]
Note that $b(y+1)=\tr(\frac{y+1}{2}(b+\sqrt{\delta}))$. We combine
the two terms involving $\psi$ to get
\[
{\small\begin{split} W_{\gL,\psi}'\left(\left(\begin{smallmatrix}
  y & 0 \\
  0 & 1 \\
\end{smallmatrix}\right)\right) = -(q-1)q^{-1}\sum_{b \in \mathbb{F}}
\sum_{z \in \mathbb{E}^\times, \atop \nr\,z =y}
\gL(z(b+\sqrt{\delta}))\psi(\tr((z+\tfrac{y+1}{2})(b+\sqrt{\delta})))+
(q-1)\delta_{y,1}.
\end{split}}
\]
Set $W_{\gL,\psi} = (q-1)^{-1}W_{\gL,\psi}'$ and $W_\gL =
\sum_{\psi \ne \psi^0} W_{\gL,\psi}$, which is a right
$K$-invariant function on $G$ belonging to $\bigoplus_{\psi \ne
\psi^0} \cl_\psi (\pi, G/K)$. We compute
\[
{\small\begin{aligned} W_\gL\left(\left(\begin{smallmatrix}
  y & 0 \\
  0 & 1 \\
\end{smallmatrix}\right)\right) &= -q^{-1}\sum_{a \in \mathbb{F}^\times} \sum_{b \in
\mathbb{F}} \sum_{z \in \mathbb{E}^\times, \atop \nr\,z =y}
\gL(z(b+\sqrt{\delta}))\psi^a(\tr((z+\tfrac{y+1}{2})(b+\sqrt{\delta})))
+ (q-1)\delta_{y,1} \\
&=-q^{-1}\sum_a \sum_b \sum_{z,\, \nr\,z = y}
\gL(za(b+\sqrt{\delta}))\psi(\tr((z+\tfrac{y+1}{2})a(b+\sqrt{\delta})))
+ (q-1)\delta_{y,1}
\end{aligned}}
\]
since $\gL$ is trivial on $\mathbb{F}^\times$. Observe that as $a$
runs through all elements in $\mathbb{F}^\times$ and $b$ runs
through all elements in $\mathbb{F}$, $a(b+\sqrt{\delta})$ runs
through all elements in $\mathbb{E} \smallsetminus \mathbb{F}$.
Therefore
\[
{\small\begin{aligned} &\sum_{a \in \mathbb{F}^\times} \sum_{b \in
\mathbb{F}} \gL(za(b+\sqrt{\delta}))
\psi(\tr((z+\tfrac{y+1}{2})a(b+\sqrt{\delta}))) \\
&=\sum_{w \in \mathbb{E}^\times}
\gL(zw)\psi(\tr((z+\tfrac{y+1}{2})w)) - \sum_{a \in
\mathbb{F}^\times} \gL(za)\psi(a\tr(z+\tfrac{y+1}{2})).
\end{aligned}}
\]
Since $\gL$ is trivial on $\mathbb{F}^\times$, the last sum is
equal to $(q-1)\gL(z)$ if $\tr(z+\frac{y+1}{2}) = \tr(z) + y +1 = 0$, and to $-\gL(z)$ otherwise. This gives
\[
{\small\begin{aligned} W_\gL\left(\left(\begin{smallmatrix}
  y & 0 \\
  0 & 1 \\
\end{smallmatrix}\right)\right) &= -q^{-1}\sum_{z \in \mbfQ, \atop
\nr\,z=y}\sum_{w \in \mbfQ^\times}
\gL(zw)\psi(\tr(z+\tfrac{y+1}{2})w) +
\sum_{z \in \mbfQ, \nr\, z =y, \atop \tr(z)+y+1=0} \gL(z) \\
&\quad +q^{-1}\sum_{z \in \mbfQ, \atop \nr\,z =y} \gL(z) + (q-1)\delta_{y,1} \\
&= -q^{-1}\sum_{z \in \mbfQ, \nr\, z =y, \atop
z+\frac{y+1}{2}\ne0} \gL(z)\bar{\gL}(z+\tfrac{y+1}{2}) \sum_{w
\in \mbfQ}\gL((z+\tfrac{y+1}{2})w)\psi(\tr((z+\tfrac{y+1}{2})w)) \\
&\quad -q^{-1} \sum_{z \in \mbfQ, \nr\,z=y \atop
z+\frac{y+1}{2}=0} \sum_{w \in \mbfQ} \gL(zw) + \sum_{z \in \mbfQ,
\nr\,z = y, \atop \tr(z)+y+1=0} \gL(z) - q^{-1}\sum_{z \in \mbfQ,
\atop \nr\,z =y} \gL(z) + (q-1)\delta_{y,1}.
\end{aligned}}
\]
Further, since $\gL$ is a nontrivial character on the kernel of
norm in $\mbfQ^\times$, we have $\sum_{\nr\,z =y} \gL(z) = 0$ and
$\sum_{w \in \mbfQ^\times} \gL(zw) =0$. Therefore
\[
{\small\begin{split} W_\gL\left(\left(\begin{smallmatrix}
  y & 0 \\
  0 & 1 \\
\end{smallmatrix}\right)\right) = - q^{-1}\gG(\gL,\psi \circ
\tr) \sum_{z \in \mbfQ, \nr\,z =y, \atop z+\frac{y+1}{2}\ne 0}
\bar{\gL}(1+\tfrac{y+1}{2z}) - \sum_{z \in \mbfQ, \nr\,z=y, \atop
\tr(z)+y+1=0} \gL(z) +(q-1)\delta_{y,1}.\end{split}}
\]
We discuss the second sum $\sum_{z \in \mbfQ, \nr\, z=y, \atop
\tr(z)+y+1=0} \gL(z)$. Since the polynomial $x^2+(y+1)x+y$ factors
as $(x+y)(x+1)$, there are no elements $z$ in $\mbfQ
\smallsetminus \mathbb{F}$ with $\nr\,z =y$ and $\tr \,z = -y-1$.
If $z \in \mathbb{F}$, then the conditions yield $y=z^2$ and
$2z=-y-1=-z^2-1$, implying the only nonvoid sum occurs when $z=-1$
and $y=1$, in which case the sum is equal to 1. As for the first
sum, in order that $z+\frac{y+1}{2}=0$ for some $z$ with
$\nr\,z=y$, we must have $z=-\frac{y+1}{2} \in \mathbb{F}^\times$.
Then $y = \nr(z) = \frac{y+1}{2}$ implies $y=1$. In this case
$z=-1$. Hence
\begin{equation}\label{Wid}
\begin{aligned}
W_\gL\left(\left(\begin{smallmatrix}
  1 & 0 \\
  0 & 1 \\
\end{smallmatrix}\right)\right) &= -q^{-1}\gG(\gL, \psi\circ\tr)
\sum_{z \in \mbfQ, z \ne -1, \atop \nr\,z=1}
\bar{\gL}(1+\tfrac{1}{z}) + 1 +(q-1) \\
&= -q^{-1}\gG(\gL, \psi\circ\tr) \sum_{z \in \mbfQ, z \ne -1,
\atop \nr\,z=1} \gL(1+z) + q,
\end{aligned}
\end{equation}
and for $y \ne 1$,
\[
W_\gL\left(\left(\begin{smallmatrix}
  y & 0 \\
  0 & 1 \\
\end{smallmatrix}\right)\right) = -q^{-1}\gG(\gL, \psi \circ \tr)
\sum_{z \in \mathbb{E}, \atop \nr\,z=y}
\bar{\gL}(1+\tfrac{y+1}{2z}) =
W_\gL\left(\left(\begin{smallmatrix}
  y^{-1} & 0 \\
  0 & 1 \\
\end{smallmatrix}\right)\right).
\]
\begin{prop}
$\sum_{\nr\,z =1 \atop z \ne -1} \gL(1+z) = -\gL(\sd) = \pm 1.$
\end{prop}
\begin{proof}
Write $\mathcal{N}$ for the subgroup of elements in
$\mathbb{E}^\times$ with norm 1 to $\mathbb{F}$. It is cyclic of
order $q+1$. The map $\phi : E^\times \to \mathcal{N}$ given by $z
\mapsto \frac{z}{z^q}$ is surjective with kernel
$\mathbb{F}^\times$. Any character of $E^\times$ trivial on
$\mathbb{F}^\times$ factors through $\mathcal{N}$. Thus there is a
character $\chi$ of $\mathcal{N}$ such that $\gL(z) =
\chi(\phi(z)) = \chi(\frac{z}{z^q})$. Consider the restriction of
$\phi$ to the subset $S = \{1+z : z \in \mathcal{N}, z \ne -1 \}$.
We claim that $\phi$ is injective on $S$. Indeed, if $z,w \in
\mathcal{N} \smallsetminus \{-1\}$ are such that $\phi(1+z) =
\phi(1+w)$, then there exists $k \in \mathbb{F}^\times$ such that
$1+z = k(1+w)$. Then $z = kw+k-1$ implies
\[
1 = z z^q = k^2w w^q + (k-1)^2 + k(k-1)(w+ w^q) = k^2 + (k-1)^2 +
k(k-1)(w + w^q).
\]
If $z \ne w$, then $k \ne 1$, and the above implies $w+w^q = -2$,
that is, $w = -1$, a contradiction. Thus $\phi(S)$ contains all
elements in $\mathcal{N}$ except $-1$, for $1+z = -(1+z^q)$ would
imply $z = -1$. Therefore
\[
\sum_{\nr \, z = 1 \atop z \ne -1} \gL(1+z) = \sum_{z \in
\mathcal{N} \atop z \ne -1} \chi(z) = -\chi(-1) = -\gL(\sd) = \pm
1
\]
since $\gL^2(\sd) = \gL(\gd) = 1$.
\end{proof}
\medskip

On the other hand,
\[
{\small\begin{split} &\gG(\gL, \psi \circ \tr) = \sum_{z \in
\mathbb{E}^\times} \gL(z) \psi(\tr\, z)
\\
&= \sum_{b \in \mathbb{F}^\times} \sum_{a \in \mathbb{F}^\times}
\gL(a(b+\sd))\psi(a\tr(b+\sd)) + \sum_{a \in \mathbb{F}^\times}
\gL(a\sd)\psi(a\tr\,\sd) + \sum_{a \in \mathbb{F}^\times}
\gL(a)\psi(a \tr\,1)
\\
&= - \sum_{b \in \mathbb{F}^\times} \gL(b+\sd) - \gL(1) + (q-1)
\gL(\sd) = q\gL(\sd).
\end{split}}
\]
Therefore $\gG(\gL, \psi \circ \tr) \sum_{\nr\,z = 1 \atop z \ne
-1} \gL(1+z) = q\gL(\sd)(-\gL(\sd)) = -q$. Plugging into
(\ref{Wid}), we obtain

\begin{prop}\label{Wat1}
$W_\Lambda \left(\left(\begin{smallmatrix}
  1 & 0 \\
  0 & 1 \\
\end{smallmatrix}\right)\right) = q+1.$
\end{prop}

Consequently $W_\gL \ne 0$, and hence $W_{\gL, \psi} \ne 0$ for
some $\psi$. We have shown that the dimension of
$\cl_\psi(\pi_\gL, G/K)$ is at least one for some and hence for
all $\psi \ne \psi^0$. We record this in

\begin{prop}
For each discrete series representation $\pi_\gL$ of $G$ and each
nontrivial additive character $\psi$ of $\mathbb{F}$, the space
$\cl_\psi(\pi_\gL, G/K)$ is $1$-dimensional.
\end{prop}

Let $KsK$ be a $K$-double coset of $G$. The operator $T_{KsK}$
preserves each space $\cl_\psi(\pi, G/K)$ and the eigenvalue
depends only on  the representation, not its model. As computed in
\cite{AC92}, the $K$-double cosets $KsK$ with cardinality greater
than $q+1$ are symmetric, and they are parameterized by $c \in
\mathbb{F}$ with $c \ne \pm 1$ so that $KsK = \bigcup_{(y,x)}
\left(\begin{smallmatrix}
  y & \delta x \\
  0 & 1 \\
\end{smallmatrix}\right)K$ where $(y,x)$ runs all solutions of
 $(y+c)^2 - \delta x^2 =
c^2-1$ over $\mathbb{F}$. Here $s$ may be any coset representative
$\left(\begin{smallmatrix}
  y & \delta x \\
  0 & 1 \\
\end{smallmatrix}\right)$, for instance. Denote this double coset by $K_c$ for short. We proceed
to compute the eigenvalues of $T_{K_c}$ using the right
$K$-invariant eigenfunctions obtained above.

For $\mu \in \mathbb{F}^\times$ with $\mu \ne 1$, we have $T_{K_c}
f_\mu = \gl_{\pi_\mu,c} f_\mu$. Hence
\[
\gl_{\pi_\mu, c} = \gl_{\pi_\mu, c} f_\mu
\left(\left(\begin{smallmatrix}
  1 & 0 \\
  0 & 1 \\
\end{smallmatrix}\right)\right) = \sum_{y \in \mathbb{F}, \atop (y+c)^2 - \gd x^2 =
c^2-1} f_\mu \left(\left(\begin{smallmatrix}
  y & \delta x \\
  0 & 1 \\
\end{smallmatrix}\right)\right) = \sum_{y \in \mathbb{F}, \atop (y+c)^2 - \gd x^2 =
c^2-1}\mu(y),
\]
which is known to have absolute value bounded by $2 \sqrt q$. See
Theorem 10 in Chapter 9 of \cite{L96} for a proof using id\`{e}le
class characters. As remarked before, the eigenvalue
$\gl_{\pi_\mu, c}$ depends only on the representation and the
double coset $K_c$. We have shown

\begin{prop} \label{muK}
Let $\mu$ be a nontrivial character of $\mathbb{F}^\times$. The
eigenvalue $\gl_{\pi_\mu, c}$ on $\cl_\psi(\pi_\mu, G/K)$ is
$\sum_{(y,x) \atop (y+c)^2-\gd x^2=c^2-1} \mu(y)$, which has
absolute value at most $2 \sqrt q$ for all additive character
$\psi$ of $\mathbb{F}$. For $\mu =1$, the trivial character, the
constant functions on $G$ are the eigenfunctions of $T_{K_c}$ with
eigenvalue $q+1$, coming from the trivial representation of $G$.
\end{prop}

Next we fix a discrete series representation $\pi_\gL$ and discuss
the eigenvalue of $T_{K_c}$ on the 1-dimensional space
$\cl_\psi(\pi_\gL, G/K)$ for $\psi \ne \psi^0$. Since
$W_{\gL,\psi} \in \cl_\psi(\pi_\gL, G/K)$, we have
$T_{K_c}W_{\gL,\psi} = \gl_{\pi_\gL,c}W_{\gL, \psi}$ and hence
$T_{K_c}W_\gL = \gl_{\pi_\gL,c} W_\gL$, where $W_\gL = \sum_{\psi
\ne \psi^0} W_{\gL, \psi}$ has the Fourier expansion
\[
W_\gL\left(\left(\begin{smallmatrix}
  1 & x \\
  0 & 1 \\
\end{smallmatrix}\right)\left(\begin{smallmatrix}
  y & 0 \\
  0 & 1 \\
\end{smallmatrix}\right)\right) = \sum_{\psi \ne \psi^0} W_{\gL,\psi} \left(\left(\begin{smallmatrix}
  y & 0 \\
  0 & 1 \\
\end{smallmatrix}\right) \right) \psi(x).
\]
Recall that $W_\gL\left(\left(\begin{smallmatrix}
  1 & 0 \\
  0 & 1 \\
\end{smallmatrix}\right)\right) = q+1 \ne 0$. We will use this value to compute the eigenvalue of
$T_{K_c}$ for $c \ne \pm 1$. By definition,
\[
{\small\begin{split}
&(T_{K_c}W_\gL)\left(\left(\begin{smallmatrix}
  1 & 0 \\
  0 & 1 \\
\end{smallmatrix}\right)\right) = \sum_{(y,x) \atop (y+c)^2 - \gd x^2 =
c^2-1} W_\gL\left(\left(\begin{smallmatrix}
  y & \gd x \\
  0 & 1 \\
\end{smallmatrix}\right)\right) = \sum_{(y,x) \atop (y+c)^2 - \gd x^2 =
c^2-1} W_\gL\left(\left(\begin{smallmatrix}
  y & 0 \\
  0 & 1 \\
\end{smallmatrix}\right)\right)\psi(\gd x) \\
&= \sum_{(y,x) \atop (y+c)^2 - \gd x^2 = c^2-1} \sum_{\psi \ne
\psi^0} - q^{-1}\sum_{b \in \mathbb{F}} \sum_{ w \in
\mathbb{E}^\times \atop \nr\, w = y(b^2 - \gd)} \gL(w)\psi(\tr \,
w) \psi(b(y+1))\psi(\gd x) + \gd_{y,1}\psi(\gd x).
\end{split}}
\]
For fixed $(y, x)$, we compute
\[
{\small\begin{split} &\sum_{\psi \ne \psi^0} - q^{-1} \sum_{b \in
\mathbb{F}} \sum_{w \in \mathbb{E}^\times \atop \nr\, w =
y(b^2-\gd)} \gL(w) \psi(\tr\,w)\psi(b(y+1))\psi(\gd x) \\
&=q^{-1} \sum_{b \in \mathbb{F}} \sum_{w \in \mathbb{E}^\times
\atop \nr\, w = y(b^2-\gd)} \gL(w) - \sum_{b \in \mathbb{F}}
\sum_{w \in \mathbb{E}^\times, \nr\, w = y(b^2 -\gd) \atop \tr\, w
= -(y+1)b-\gd x}  \gL(w) = - \sum_{b \in \mathbb{F}} \sum_{w \in
\mathbb{E}^\times, \nr\, w = y(b^2 -\gd) \atop \tr\, w =
-(y+1)b-\gd x} \gL(w)
\end{split}}
\]
since $\sum_{\nr\,w=y(b^2-\gd)} \gL(w) = 0$ for all $b \in
\mathbb{F}$. Further,
\[
\sum_{\psi \ne \psi^0} \gd_{y,1} \psi(\gd x) = - \gd_{y,1} + q
\gd_{y,1} \gd_{x,0} = - \gd_{y,1}
\]
as $(y,x) = (1,0)$ does not satisfy the equation $(y+c)^2 - \gd
x^2 = c^2 -1$ because $c \ne -1$. Therefore we may write
\[
(T_{K_c}W_\gL)\left(\left(\begin{smallmatrix}
  1 & 0 \\
  0 & 1 \\
\end{smallmatrix}\right)\right) = \sum_{(y,x) \atop (y+c)^2 -
\gd x^2 = c^2 -1} S_{(y,x)},
\]
where
\begin{equation}\label{syx}
S_{(y,x)} = -\sum_{b \in \mathbb{F}} \sum_{w \in
\mathbb{E}^\times, \nr\,w=y(b^2-\gd) \atop \tr\,w=-(y+1)b -\gd x}
\gL(w) - \gd_{y,1}.
\end{equation}

To proceed, we prove

\begin{thm} \label{mainK}
Given $c \in \mathbb{F}$, $c \ne \pm 1$, let $y,x \in \mathbb{F}$
satisfy $(y+c)^2 - \gd x^2 = c^2 -1$. Then $|S_{(y,x)}| \le 2\sqrt
q$.
\end{thm}

The following two results, adapted from Theorem 4 in Chapter 6 of
\cite{L96} and Theorem 5 of \cite{L99}, respectively, will be used
repeatedly in the proof.

\begin{thm} \label{nu}
Let $\nu$ be the quadratic character of $\mathbb{F}^\times$. Let
$f(T)$ be a quadratic polynomial over $\mathbb{F}$ with two
distinct roots. Denote by $v_1, \dots, v_r, r \le 2$, the places
of $\mathbb{F}(T)$ containing the roots of $f$. Then there exists
an id\`{e}le class character $\eta = \eta_{\nu,f}$ of
$\mathbb{F}(T)$ such that
\begin{itemize}
    \item[(1)] The conductor of $\eta$ is $v_1+ \cdots + v_r;$
    \item[(2)] At each place $v$ of degree $1$ with uniformizer $\pi_v =
    T-v$ where $\eta$ is unramified,  we have $\eta_v(\pi_v) =
    \nu(f(v))$, and $\eta_\nf(\pi_\nf) = 1$.
\end{itemize}
\end{thm}

\begin{thm} \label{omega}
Let $h(T)$ be a nonconstant polynomial over $\mathbb{E}$ with
distinct roots. Let $w_1, \dots, w_r$ be the places of
$\mathbb{F}(T)$ containing the roots of $h$. Then there exists an
id\`{e}le class character $\omega = \omega_{\gL,h}$ of
$\mathbb{F}(T)$ such that
\begin{itemize}
    \item[(1)] The conductor of $\omega$ is $w_1 + \cdots + w_r;$
    \item[(2)] At each place $v$ of degree $1$ with uniformizer $\pi_v =
    T-v$ where $\omega$ is unramified, we have $\omega_v(\pi_v) =
    \gL(h(v))$.
\end{itemize}
\end{thm}

The character $\omega$ in Theorem \ref{omega} is unramified at
$\nf$ since $\gL$ is trivial on $\mathbb{F}^\times$. But the value
of $\omega_\nf (\pi_\nf)$ depends on $\gL$ and $h$. In case $h(0)
\ne 0$, it is equal to
\begin{equation}\label{oinfinity}
\begin{split}
\omega_\nf(\pi_\nf) &= \prod_{v \ne \nf} \omega_v(T) =
\gL(h(0))\gL(\text{product of roots of $h$})^{-1} \\
&= \gL(\text{the leading coefficient of $h$}),
\end{split}
\end{equation}
following the proof of Theorem 5 in \cite{L99}.

The following character sum estimate results from the Riemann
hypothesis for curves, as explained in Section 1, Chapter 6 of
\cite{L96}. It will be used repeatedly to derive character sum
estimates.

\begin{prop}
Let $\chi$ be an id\`{e}le class character of $\mathbb{F}(T)$ such
that its conductor has degree $m$. Then
\[
\left| \sum_{\deg v = 1 \atop \chi_v \, \mathrm{unramified}}
\chi_v(\pi_v) \right| \le (m-2) \sqrt q.
\]
\end{prop}

We now begin the proof of Theorem \ref{mainK}. We distinguish
three cases.

\noindent {\bf Case 1.} $y=-1$. Then $\gd x^2 = 2-2c$ implies $x
\ne 0$ since $c \ne 1$. In this case
\[
{\small\begin{split} S_{(-1,x)} = -\sum_{b \in \mathbb{F}} \sum_{w
\in \mathbb{E}^\times, \nr\,w=-(b^2-\gd) \atop \tr\,w=-\gd x}
\gL(w).\end{split}}
\]
Write $w = -\frac{\gd x}{2}+v \sd$. Then $\nr\, w =
\frac{\gd^2x^2}{4}-v^2\gd = \gd-b^2$ amounts to $\gd+v^2\gd
-\frac{\gd^2 x^2}{4} = b^2$ being a square, or equivalently,
$1+v^2-\frac{\gd x^2}{4} = 1+v^2-\frac{1-c}{2} =
v^2+\frac{1+c}{2}$ not a square in $\mathbb{F}$. Denote by $\nu$
the quadratic character of $\mathbb{F}^\times$, extended to a
function on $\mathbb{F}$ by letting $\nu(0) = 0$. We rewrite
$S_{(-1,x)}$ as
\[
{\small\begin{split} &S_{(-1,x)} = -\sum_{v \in \mathbb{F} \atop
v^2 + \frac{1+c}{2} \, \text{is a nonsquare}} \gL\left(-\tfrac{\gd x}{2}+v\sd\right) \\
&= -\frac{1}{2}\sum_{v \in \mathbb{F}} \left( 1 -\nu\left( v^2+
\tfrac{1+c}{2}\right)\right) \gL\left(- \tfrac{\gd x}{2} + v\sd
\right)+ \frac{1}{2}\sum_{v \in \mathbb{F} \atop v^2 +
\frac{1+c}{2}=0}\gL\left( -\tfrac{\gd x}{2} + v\sd \right) \\
&= \frac{1}{2}\sum_{v \in \mathbb{F}} \nu\left( v^2 +
\tfrac{1+c}{2}\right) \gL\left( -\tfrac{\gd x}{2} + v\sd \right)-
\frac{1}{2}\sum_{v \in \mathbb{F}} \gL\left(- \tfrac{\gd x}{2} +
v\sd \right) + \frac{1}{2}\sum_{v \in \mathbb{F} \atop v^2 +
\frac{1+c}{2}=0} \gL\left(- \tfrac{\gd x}{2} + v\sd \right).
\end{split}}
\]
As $v$ runs through all elements in $\mathbb{F}$, no two elements
of the form $-\frac{\gd x}{2} + v \sd$ differ by a multiple in
$\mathbb{F}^\times$, hence $-\frac{1}{2} \sum_{v \in \mathbb{F}}
\gL(-\frac{\gd x}{2} + v \sd) = \frac{1}{2}\gL(\sd)$ and
\[
{\small\begin{split}S_{(-1,x)} = \frac{1}{2}\sum_{v \in
\mathbb{F}} \nu\left( v^2 + \tfrac{1+c}{2}\right) \gL\left(
-\tfrac{\gd x}{2} + v\sd \right) + \frac{1}{2}\gL(\sd)+
\frac{1}{2}\sum_{v \in \mathbb{F} \atop v^2 + \frac{1+c}{2}=0}
\gL\left(- \tfrac{\gd x}{2} + v\sd \right).\end{split}}
\]
Let $f(T) = T^2 +\frac{1+c}{2}$ and $h(T) = -\frac{\gd x}{2} +
T\sd$. Let $\eta = \eta_{\nu, f}$ and $\omega = \omega_{\gL,h}$ be
the id\`{e}le class characters of $\mathbb{F}(T)$ as described in
Theorems \ref{nu} and \ref{omega}, respectively. The conductor of
$\omega$ is $w_1$ with the uniformizer $\pi_{w_1} =
T^2-\frac{1-c}{2}$, which is disjoint from the conductor of
$\eta$. Hence the conductor of $\eta \omega$ has degree $4$.
Moreover, $\omega_\nf(\pi_\nf) = \gL(\sd)$ by equation
(\ref{oinfinity}). Therefore
\[
{\small\begin{split} \left| \sum_{\deg v = 1 \atop \eta_v \omega_v
\, \text{unramified}} \eta_v(\pi_v) \omega_v(\pi_v) \right| &=
\left|
\sum_{v \in \mathbb{F}} \nu(f(v))\gL(h(v)) + \gL(\sd) \right| \\
&= \left| \sum_{v \in \mathbb{F}} \nu\left( v^2 +
\tfrac{1+c}{2}\right)
\gL\left(-\tfrac{\gd x}{2}+ v\sd \right) + \gL(\sd) \right| \\
&=\left| 2S_{(-1,x)} - \sum_{v \in \mathbb{F} \atop v^2 +
\frac{1+c}{2} = 0} \gL\left( -\tfrac{\gd x}{2} + v
\sd\right)\right| \le (4-2)\sqrt q
\end{split}}
\]
implies
\[
|S_{(-1,x)}| \le \sqrt q +1 < 2\sqrt q,
\]
as desired.

It remains to deal with the case $y \ne -1$. We eliminate the
variable $b$ in the expression of $S_{(y,x)}$. Write $w = u+v\sd$
with $u,v \in \mathbb{F}$. The conditions $\nr \, w = y(b^2-\gd)$
and $\tr\, w = -(y+1)b - \gd x$ can be combined as
\[
\frac{\nr\, w}{y} + \gd = b^2 = \left( \frac{\tr\, w + \gd
x}{y+1}\right)^2.
\]
In other words,
\[
u^2-v^2\gd +y\gd = \frac{y}{(y+1)^2}(2u+\gd x)^2 =
\frac{4y}{(y+1)^2}u^2 + \frac{4yu\gd x}{(y+1)^2} + \frac{y \gd^2
x^2}{(y+1)^2},
\]
which in turn yields
\begin{equation}\label{yne-1}
\left(\frac{y-1}{y+1}\right)^2u^2 - \frac{4y\gd x}{(y+1)^2}u -
v^2\gd = -y \gd + \frac{y\gd^2 x^2}{(y+1)^2}.
\end{equation}

\noindent{\bf Case 2.} $y=1$. Then $\gd x^2 = 2+2c \ne 0$ since $c
\ne -1$. The above relation simplifies as
\[
{\small\begin{split} -\gd x u - v^2\gd = -\gd +\frac{\gd^2x^2}{4}
= -\gd +\gd\frac{1+c}{2} = \gd\frac{c-1}{2},\end{split}}
\]
which allows us to express $u$ in terms of $v$:  $u = \tfrac{1}{x}
\left(\tfrac{1-c}{2} - v^2\right)$. Then
\[
{\small\begin{split} S_{(1,x)} = - \sum_{v \in \mathbb{F}}
\gL\left( \tfrac{1}{x}\left( \tfrac{1-c}{2}-v^2\right) +
v\sd\right) -1. \end{split}}
\]
Let
\[
{\small\begin{split} h(T) &= \tfrac{1}{x}\left( \tfrac{1-c}{2} -
T^2\right) +
T\sd = -\tfrac{1}{x} \left( T^2 - x\sd T -\tfrac{1-c}{2}\right) \\
&= -\tfrac{1}{x}\left(\left( T - \tfrac{x \sd}{2}\right)^2 -
\tfrac{x^2\delta}{4} + \tfrac{c-1}{2}\right) =
-\tfrac{1}{x}\left(\left( T - \tfrac{x \sd}{2}\right)^2 -
1\right).
\end{split}}
\]
Let $\go = \go_{\gL,h}$ be the id\`{e}le class character of
$\mathbb{F}(T)$ attached to $\gL$ and $h$ as described in Theorem
\ref{omega}. The conductor of $\go$ is $w_1 + w_2$, where $w_1$
and $w_2$ are two degree two places of $\mathbb{F}(T)$ containing
$\frac{x\sd}{2} - 1$ and $\frac{x \sd}{2} + 1$ as roots,
respectively. So the conductor of $\go$ has degree 4. By equation
(\ref{oinfinity}), $\go_\nf(\pi_\nf) =1$. Put together, we have
\[
-S_{(1,x)} = \sum_{v \in \mathbb{F}} \gL(h(v)) + 1 = \sum_{\deg v
=1} \go_v(\pi_v),
\]
which satisfies $|S_{(1,x)}| \le 2\sqrt q$.

\noindent{\bf Case 3.} $y \ne \pm 1.$ We have $y^2+2yc+1 = \gd
x^2$. The relation (\ref{yne-1}) can be rewritten as
\[
{\small\begin{split} \left( \frac{y-1}{y+1}u - \frac{2y\gd
x}{y^2-1} \right)^2 - v^2\gd = \gd y^2 \frac{2+2c}{(y-1)^2}
\end{split}}
\]
so that
\[
{\small\begin{split} S_{(y,x)} = - \sum_{\nr (\frac{y-1}{y+1}u -
\frac{2y\gd x}{y^2-1} + v \sd) = \gd y^2\frac{2+2c}{(y-1)^2}}
\gL(u+v\sd).\end{split}}
\]
Replacing $u$ by $\frac{y+1}{y-1}u + \frac{2y\delta x}{(y-1)^2}$,
we rewrite the above as
\[
{\small\begin{split} S_{(y,x)} = - \sum_{\nr(u+v\sd)= \gd
y^2\frac{2+2c}{(y-1)^2}} \gL\left( \frac{y+1}{y-1}u + \frac{2y\gd
x}{(y-1)^2} + v\sd\right).\end{split}}
\]
Set $z = u+v\sd$. Then
\[
{\small\begin{split} 2u = z + z^q = z + \frac{\nr\, z}{z} = z +
\frac{1}{z}\gd y^2\frac{2+2c}{(y-1)^2}.\end{split}}
\]
Using this, we may express the argument of $\gL$ as a rational
function in $z$:
\[
{\small \begin{split} \frac{y+1}{y-1} u + \frac{2y\gd x}{(y-1)^2}
+ v\sd &= z+ u\frac{2}{y-1} + \frac{2y\gd x}{(y-1)^2} \\
&= \frac{y}{y-1}z + \frac{1}{z}\frac{2+2c}{(y-1)^3}\gd y^2 +
\frac{2y\gd x}{(y-1)^2} =: Q(z).
\end{split}}
\]
By choosing an element $w \in \mathbb{E} \smallsetminus
\mathbb{F}$ with $\nr\, w = \gd y^2 \frac{2+2c}{(y-1)^2}$, we
rewrite $S_{(y,x)}$ as
\[
{\small\begin{split} S_{(y,x)} = - \sum_{\nr\, z = \gd
y^2\frac{2+2c}{(y-1)^2}} \gL(Q(z)) = - \sum_{\nr\, z_1 = 1}
\gL(Q(wz_1)).\end{split}}
\]
Here
\[
{\small\begin{split} Q(wz_1) = \frac{y}{y-1} wz_1 +
\frac{1}{z_1w}\frac{\nr\,w}{y-1}\gd y^2 + \frac{2y\gd x}{(y-1)^2}
 = \frac{1}{y-1}\left(ywz_1 + \frac{w^q}{z_1} + \frac{2y\gd
x}{y-1}\right) =: R(z_1).\end{split}}
\]
Consider
\[
{\small \begin{split} &R\left(\frac{T-\sd}{T+\sd}\right) =
\frac{1}{y-1}\left(yw\frac{T-\sd}{T+\sd} +
\frac{w^q(T+\sd)}{T-\sd} + \frac{2y \gd x}{y-1}\right) \\
&= \frac{1}{(y-1)(T^2-\gd)}\left( yw(T-\gd)^2 + w^q(T+\sd)^2 +
\frac{2y\gd x}{y-1}(T^2-\gd)\right) =:
\frac{h(T)}{(y-1)(T^2-\gd)},
\end{split}}
\]
where $h(T) = ( yw + w^q + \frac{2y\gd x}{y-1})T^2 + (2w^q\sd -
2yw\sd)T + yw\gd + w^q\gd - \frac{2y\gd^2x}{y-1} \in
\mathbb{E}[T]$. As $T$ runs through elements in $\mathbb{F}$,
$\frac{T-\sd}{T+\sd}$ runs through all elements $z_1$ in $E$ with
norm 1 except $z_1 = 1$. Observe that for $T \in \mathbb{F}$,
$\gL(R(\frac{T-\sd}{T+\sd}))= \gL(\frac{h(T)}{(y-1)(T^2-\gd)}) =
\gL(h(T))$ since $\gL$ is trivial on $\mathbb{F}^\times$. Thus we
have
\[
{\small\begin{split} S_{(y,x)} = -\sum_{\nr\,z_1=1} \gL(R(z_1)) =
- \sum_{v \in \mathbb{F}} \gL(h(v)) - \gL(R(1)).\end{split}}
\]
Notice that in the course of deriving various expressions of
$S_{(y,x)}$, the character $\gL$ is always evaluated at nonzero
elements in $\mathbb{E}$. In particular, this means that $h(v) \ne
0$ for all $v \in \mathbb{F}$. In other words, the roots of $h$
are outside $\mathbb{F}$. Since $y \ne 1$ and $w \notin
\mathbb{F}$ by choice, both the leading coefficient and the
constant term of $h(T)$ lie in $\mathbb{E} \smallsetminus
\mathbb{F}$. Moreover, their ratio also lies in $\mathbb{E}
\smallsetminus \mathbb{F}$. This implies that if $h(T)$ has a root
$t$ of multiplicity 2, then $t$ lies in $\mathbb{E}\smallsetminus
\mathbb{F}$. In this case, for $v \in \mathbb{F}$,
\[
{\small\begin{split} \gL(h(v)) = \gL\left(yw + w^q + \frac{2y\gd
x}{y-1} \right)\gL^2(v-t).\end{split}}
\]
As $v$ runs through all elements in $\mathbb{F}$, no two elements
of the form $v-t$ differ by a multiple in $\mathbb{F}^\times$.
Since the order of $\gL$ is greater than 2, we have
\[
{\small \begin{split} -\sum_{v \in \mathbb{F}} \gL(h(v)) =
\gL\left(yw + w^q + \frac{2y\gd x}{y-1}\right)\gL^2(1)
\end{split}}
\]
and consequently $|S_{(y,x)}| \le 2 < 2\sqrt q$.

Finally we discuss the case where $h$ has two distinct roots.
Either they are nonconjugate over $\mathbb{F}$, and hence
contained in two distinct places $v_1, v_2$ of $\mathbb{F}(T)$ of
degree 2, or they are conjugate over $\mathbb{F}$ and contained in
a degree 4 place $v_1$ of $\mathbb{F}(T)$. At any rate, the
id\`{e}le class character $\go = \go_{\gL,h}$ attached to $\gL$
and $h$ as described in Theorem \ref{omega} has conductor of
degree 4. Further, $\go$ is unramified at all places of degree 1,
and
\[
{\small \begin{split}\go_\nf(\pi_\nf) = \gL\left(yw + w^q +
\frac{2y\gd x}{y-1}\right) = \gL(R(1)).\end{split}}
\]
Therefore
\[
{\small \begin{split} S_{(y,x)} = -\sum_{v \in \mathbb{F}}
\gL(h(v)) - \gL(R(1)) = - \sum_{\deg v =1} \go_v(\pi_v),
\end{split}}
\]
and $|S_{(y,x)}| \le 2 \sqrt q$. This completes the proof of the
theorem.

\bigskip

Since
\[
{\small\begin{split} (T_{K_c}W_\gL)\left(\left(\begin{smallmatrix}
  1 & 0 \\
  0 & 1 \\
\end{smallmatrix}\right)\right) = \sum_{(y,x) \atop (y+c)^2 - \gd x^2 = c^2
-1}S_{(y,x)} = \gl_{\pi_\gL,c}W_\gL\left(\left(\begin{smallmatrix}
  1 & 0 \\
  0 & 1 \\
\end{smallmatrix}\right)\right) = \gl_{\pi_\gL, c}(q+1)
\end{split}}
\]
is a sum of $q+1$ terms $S_{(y,x)}$ and each $S_{(y,x)}$ is of
absolute value at most $2\sqrt q$, we conclude that the eigenvalue
$\gl_{\pi_\gL,c}$ of $T_{K_c}$ on the space $\cl_\psi(\pi_\gL,
G/K)$ associated to the discrete series $\pi_\gL$ satisfies
$|\gl_{\pi_\gL,c}| \le 2 \sqrt q$ for each nontrivial additive
character $\psi$. This proves

\begin{prop}
Let $\gL$ be a character of $\mathbb{E}^\times$ trivial on
$\mathbb{F}^\times$ and of order greater than $2$. The eigenvalue
$\gl_{\pi_\gL, c}$ of the operator $T_{K_c}$ on the
$1$-dimensional space $\cl_\psi(\pi_\gl, G/K)$ is equal to
$\frac1{q+1} \sum_{(y,x) \atop (y+c)^2 - \gd x^2 = c^2
-1}S_{(y,x)}$ and satisfies $|\gl_{\pi_\gL,c}| \le 2\sqrt q$ for
all nontrivial additive character $\psi$.
\end{prop}

Combined with proposition \ref{muK}, we obtain another proof of
the following result established by Terras et al in \cite{AC92}
and \cite{CP93}.

\begin{thm}
Given $c \in \mathbb{F}$ with $c \ne \pm 1$, denote by $K_c$ the
$K$-double coset $K \left(\begin{smallmatrix}
y & \delta x \\
0 & 1 \\
\end{smallmatrix}\right)K = \cup_{(y,x)} \left(\begin{smallmatrix}
y & \delta x \\
0 & 1 \\
\end{smallmatrix}\right)K$ where $(y,x)$ satisfies $(y+c)^2 - \gd x^2 = c^2 -1$.
Then the Cayley graph $X_c = \mathrm{Cay}(G/K, K_c)$ is an
undirected $(q+1)$-regular Ramanujan graph.
\end{thm}

\bigskip

\section{The $U$-graphs}

In this section, we let $H=U$. We first analyze the space
$\cl_{\psi^0}(G/U)$. Note that $G = AU \cup
AU\left(\begin{smallmatrix}
0 & 1 \\
-1 & 0 \\
\end{smallmatrix}\right)U$ is a disjoint union of two double
cosets. One sees immediately that for each $\mu \in \hf^\times$,
the right $U$-invariant subspace of $\Ind\mu$ is 2-dimensional,
generated by
\[
g_\mu\left(\left(\begin{smallmatrix}
  y & 0 \\
  0 & 1 \\
\end{smallmatrix}\right)U\right) = \mu(y) \; \text{for all $y \in
\mathbb{F}^\times$ and} \;  g_\mu
\left(AU\left(\begin{smallmatrix}
0 & 1 \\
-1 & 0 \\
\end{smallmatrix}\right)U\right) =0.
\]
and
\[
h_\mu\left(\left(\begin{smallmatrix}
  y & x \\
  0 & 1 \\
\end{smallmatrix}\right)\left(\begin{smallmatrix}
  0 & 1 \\
  -1 & 0 \\
\end{smallmatrix}\right)U\right) = \mu(y) \; \text{for all $y \in
\mathbb{F}^\times, x \in \mathbb{F}$,  and} \;  h_\mu (AU) =0.
\]
Thus $\cl_{\psi^0}(G/U)$ is $2(q-1)$-dimensional. Next we fix a
nontrivial additive character $\psi$ of $\mathbb{F}$. For each
irreducible representation $\pi$ of $G$ with $\deg \pi >1$ and
each $v \in \ck_\psi(\pi)$, the Whittaker function
\[
\bW_v(g) = q^{-1}\sum_{s \in \mathbb{F}} (\pi(gu_s)v)(1)
\]
is right $U$-invariant. Here $u_s = \left(\begin{smallmatrix}
  1 & s \\
  0 & 1 \\
\end{smallmatrix}\right)$ as in \S 2. Using the actions of
$\pi(u_s)$ on $\gh \in \hf^\times, D_0$, and $D_\infty$ described
in \S 2, one gets
\[
\bW_\gh(g) = 0, \bW_{D_0}(g) = (\pi(g)D_0)(1) \; \text{and} \;
\bW_{D_\infty}(g) = (\pi(g)D_\infty)(1)
\]
for all $g \in G$. Since $\dim \cl(G/U) = (q+1)(q-1)$, by
dimension counting one finds that the right $U$-invariant space is
$0$-dimensional for discrete series representations; it is
$1$-dimensional generated by $\{\bW_{D_0}\}$ for the Steinberg
representations $\pi_\mu$ with $\mu^2=1$, and it is
$2$-dimensional generated by $\{\bW_{D_0},\bW_{D_\infty}\}$ for
principal series representations $\pi_\mu$ with $\mu^2 \ne 1$. We
record this in

\begin{prop}
\begin{itemize}
    \item[(1)] The space $\cl_{\psi^0}(\pi_\mu,G/U)$ is
    $2$-dimensional for all $\mu \in \hf^\times$.
    \item[(2)] For nontrivial $\psi$'s, the space $\cl_\psi(\pi,
    G/U)$ is $0$-dimensional if $\pi$ is a discrete series representation$;$ it is
    $1$-dimensional if $\pi$ is a Steinberg representation, and $2$-dimensional
    if $\pi$ is a principal series representation
\end{itemize}
\end{prop}

The $2(q-1)$ $U$-double cosets of $G$ fall in two categories. The
first consists of $U\left(\begin{smallmatrix}
  r & 0 \\
  0 & 1 \\
\end{smallmatrix}\right)U, r \in \mathbb{F}^\times$, which are not symmetric if $r \ne \pm1$
and which are contained in the Borel subgroup $(=AU)$ of $G$,
hence not interesting; the second consists of
$U\left(\begin{smallmatrix}
  0 & t \\
  -1 & 0 \\
\end{smallmatrix}\right)U, t\in \mathbb{F}^\times$, which are symmetric and of interest
to us. Write $U_t = U\left(\begin{smallmatrix}
  0 & t \\
  -1 & 0 \\
\end{smallmatrix}\right)U$ for short. Then
\[
U_t = \bigcup_{c \in \mathbb{F}} \left(\begin{smallmatrix}
c & t \\
-1 & 0 \\
\end{smallmatrix}\right)U
\]
is a disjoint union of $q$ $U$-cosets. The operator $T_{U_t}$
preserves each space $\cl_\psi(\pi, G/U)$. To study its
eigenvalues and eigenfunctions, we start with $\psi^0$. Note that
\[
\left(\begin{smallmatrix}
c & t \\
-1 & 0 \\
\end{smallmatrix}\right)U = \left(\begin{smallmatrix}
 t & -c \\
 0 & 1 \\
\end{smallmatrix}\right)\left(\begin{smallmatrix}
0 & 1 \\
-1 & 0 \\
\end{smallmatrix}\right)U
\]
for all $c \in \mathbb{F}$ and
\[
\left(\begin{smallmatrix}
0 & 1 \\
-1 & 0 \\
\end{smallmatrix}\right)\left(\begin{smallmatrix}
0 & t \\
-1 & 0 \\
\end{smallmatrix}\right)U = \left(\begin{smallmatrix}
t^{-1} & 0 \\
0 & 1 \\
\end{smallmatrix}\right)U \; \text{and} \; \left(\begin{smallmatrix}
0 & 1 \\
-1 & 0 \\
\end{smallmatrix}\right)\left(\begin{smallmatrix}
c & t \\
-1 & 0 \\
\end{smallmatrix}\right)U = \left(\begin{smallmatrix}
tc^{-2} & t c^{-1} \\
0 & 1 \\
\end{smallmatrix}\right)\left(\begin{smallmatrix}
0 & 1 \\
-1 & 0 \\
\end{smallmatrix}\right)U
\]
for all $c \in \mathbb{F}^\times$. Let $\mu \in \hf^\times$. Then
for $l_\mu \in \{g_\mu,h_\mu\}$, we have
\[
{\small\begin{split}(T_{U_t}l_\mu)(U) = \sum_{c\in \mathbb{F}}
l_\mu \left(\left(\begin{smallmatrix}
c & t \\
-1 & 0 \\
\end{smallmatrix}\right)U\right) = \sum_{c \in \mathbb{F}} l_\mu \left(\left(\begin{smallmatrix}
t & -c \\
0 & 1 \\
\end{smallmatrix}\right)\left(\begin{smallmatrix}
0 & 1 \\
-1 & 0 \\
\end{smallmatrix}\right)U\right)
= q \mu(t) l_\mu\left(\left(\begin{smallmatrix}
0 & 1 \\
-1 & 0 \\
\end{smallmatrix}\right)U\right)\end{split}}
\]
and
\[
{\small\begin{aligned}
(T_{U_t}l_\mu)\left(\left(\begin{smallmatrix}
0 & 1 \\
-1 & 0 \\
\end{smallmatrix}\right)U\right) &= \sum_{c \in \mathbb{F}} l_\mu \left(\left(\begin{smallmatrix}
0 & 1 \\
-1 & 0 \\
\end{smallmatrix}\right)\left(\begin{smallmatrix}
c & t \\
-1 & 0 \\
\end{smallmatrix}\right)U\right)
=  l_\mu\left(\left(\begin{smallmatrix}
t^{-1} & 0 \\
0 & 1 \\
\end{smallmatrix}\right)U\right) +  \sum_{c \in \mathbb{F}^\times} l_\mu \left(\left(\begin{smallmatrix}
tc^{-2} & c^{-1} \\
0 & 1 \\
\end{smallmatrix}\right)\left(\begin{smallmatrix}
0 & 1 \\
-1 & 0 \\
\end{smallmatrix}\right)U\right) \\
&= \mu(t^{-1})l_\mu(U) + \mu(t) \sum_{c \in \mathbb{F}^\times}
\mu^2(c) l_\mu \left(\left(\begin{smallmatrix}
0 & 1 \\
-1 & 0 \\
\end{smallmatrix}\right)U\right).
\end{aligned}}
\]
Thus
\[
T_{U_t} g_\mu = \mu(t^{-1})h_\mu \quad \text{and} \quad T_{U_t}
h_\mu = q \mu(t) g_\mu + \mu(t) \sum_{c \in \mathbb{F}^\times}
\mu^2(c) h_\mu.
\]
In other words, with respect to the basis $\{g_\mu, h_\mu\}$ of
$\cl_{\psi^0}(\pi_\mu,G/U)$, the operator $T_{U_t}$ can be
represented by the matrix
\[
{\footnotesize \left(
\begin{matrix}
  0 & q\mu(t) \\
  \mu(t^{-1}) & ~\mu(t)\sum_{c\in \mathbb{F}^\times} \mu^2(c) \\
\end{matrix}
\right).}
\]
Consequently, if $\mu^2 \ne 1$, the eigenvalues of $T_{U_t}$ are
$\pm \sqrt{q}$ with eigenfunctions $\pm \sqrt{q}\mu(t) g_\mu +
h_\mu$; while if $\mu^2 = 1$, the eigenvalues are $\mu(t)q$ and
$-\mu(t)$ with eigenfunctions $g_\mu + h_\mu$ and $q g_\mu -
h_\mu$, respectively.

Next we deal with the case $\psi \ne \psi^0$. Observe that
\[
\left(\begin{smallmatrix}
  c & t \\
  -1 & 0 \\
\end{smallmatrix}\right) = \left(\begin{smallmatrix}
  1 & -c \\
  0 & 1 \\
\end{smallmatrix}\right)\left(\begin{smallmatrix}
  0 & 1 \\
  -1 & 0 \\
\end{smallmatrix}\right)\left(\begin{smallmatrix}
  t^{-1} & 0 \\
  0 & 1 \\
\end{smallmatrix}\right)\left(\begin{smallmatrix}
  t & 0 \\
  0 & t \\
\end{smallmatrix}\right)
\]
and $\gG(1, \psi) = -1$. Recall that for a Steinberg
representation $\pi= \pi_\mu$ with $\mu^2 =1$, the space
$\cl_\psi(\pi, G/U)$ is 1-dimensional generated by
$\{\bW_{D_0}\}$. The action of $T_{U_t}$ on $\bW_{D_0}$ is given
by, according to the Kirillov model of $\pi_\mu$,
{\small\begin{align} &(T_{U_t}\bW_{D_0})(g) = \sum_{c \in
\mathbb{F}} \bW_{D_0} \left(g\left(\begin{smallmatrix}
  c & t \\
  -1 & 0 \\
\end{smallmatrix}\right)\right) = \sum_{c \in \mathbb{F}}\left(\pi\left(g\left(\begin{smallmatrix}
  c & t \\
  -1 & 0 \\
\end{smallmatrix}\right)\right)D_0\right)(1) \notag \\
&= \sum_{c \in \mathbb{F}}
(\pi(g)\pi(u_c)\pi(w)\pi(h_{t^{-1}})D_0)(1) = \sum_{c \in
\mathbb{F}} (\pi(g)\pi(u_c)\pi(w)\mu(t)D_0)(1) \quad
\text{(since $\mu = \mu^{-1}$)} \notag \\
&= \sum_{c \in \mathbb{F}} (\pi(g)\pi(u_c)(-q^{-1} \ve(\pi,
\mu^{-1},
\psi)\mu(t))(\mu^{-1}+D_0))(1) \notag \\
&= -q^{-1} \ve(\pi, \mu^{-1}, \psi)\mu(t) \pi(g)\left[ \sum_{c\in
\mathbb{F}^\times} (q-1)^{-1} \sum_{\beta \in \hf^\times} \beta
\mu(c) \gG(\mu^{-1}\beta^{-1}, \psi)\beta + \mu^{-1} + \sum_{c \in
\mathbb{F}} D_0
\right](1) \notag \\
&= -q^{-1} \ve(\pi, \mu^{-1}, \psi)\mu(t) \pi(g)\left[ (q-1)^{-1}
(q-1)\gG(1, \psi)\mu^{-1} + \mu^{-1} + qD_0\right](1) \notag\\
&= -\ve(\pi, \mu^{-1}, \psi) \mu(t) [\pi(g)D_0](1) = -\ve(\pi,
\mu^{-1}, \psi) \mu(t) \bW_{D_0}(g) \notag
\end{align}}
for all $g \in G$, which implies that the eigenvalue is
\[
 - \ve(\pi, \mu^{-1}, \psi)\mu(t) = - \gG(\mu \mu^{-1},
\psi)\gG(\mu^{-1} \mu^{-1}, \psi) \mu(t) = - \gG(1, \psi)\gG(1,
\psi) \mu(t) = -\mu(t).
\]
Recall that for a principal series representation $\pi = \pi_\mu$
with $\mu^2 \ne 1$, the space $\cl_\psi(\pi_\mu, G/U)$ is
$2$-dimensional generated by $\{\bW_{D_0}, \bW_{D_\infty}\}$.
Using the Kirillov model of $\pi$, we get

{\small\begin{align} &(T_{U_t}\bW_{D_0})(g) = \sum_{c \in
\mathbb{F}} \bW_{D_0} \left(g\left(\begin{smallmatrix}
  c & t \\
  -1 & 0 \\
\end{smallmatrix}\right)\right) = \sum_{c \in \mathbb{F}}\left(\pi\left(g\left(\begin{smallmatrix}
  c & t \\
  -1 & 0 \\
\end{smallmatrix}\right)\right)D_0\right)(1)\notag \\
&= \sum_{c \in \mathbb{F}}
(\pi(g)\pi(u_c)\pi(w)\pi(h_{t^{-1}})D_0)(1)=
\sum_{c \in \mathbb{F}} (\pi(g)\pi(u_c)\pi(w)\mu(t^{-1})D_0)(1)\notag\\
&= \sum_{c \in \mathbb{F}} (\pi(g)\pi(u_c)(-q^{-1} \ve(\pi, \mu,
\psi)\mu(t^{-1}))(\mu^{-1}+D_0))(1)\notag \\
&= -q^{-1} \ve(\pi, \mu, \psi)\mu(t^{-1}) \pi(g)\left[ \sum_{c\in
\mathbb{F}^\times} (q-1)^{-1} \sum_{\beta \in \hf^\times} \beta
\mu(c) \gG(\mu^{-1}\beta^{-1}, \psi)\beta + \mu^{-1} + \sum_{c \in
\mathbb{F}}
D_\infty\right](1) \notag\\
&= -q^{-1} \ve(\pi, \mu, \psi)\mu(t^{-1}) \pi(g)\left[ (q-1)^{-1}
(q-1)\gG(1, \psi)\mu^{-1} + \mu^{-1} + qD_\infty\right](1) \notag\\
&= -\ve(\pi, \mu, \psi) \mu(t^{-1}) [\pi(g)D_\infty](1) =
-\ve(\pi, \mu, \psi) \mu(t^{-1}) \bW_{D_\infty}(g)\notag
\end{align}}
and {\small\begin{align} &(T_{U_t}\bW_{D_\infty})(g) = \sum_{c \in
\mathbb{F}} \bW_{D_\infty} \left(g\left(\begin{smallmatrix}
  c & t \\
  -1 & 0 \\
\end{smallmatrix}\right)\right) = \sum_{c \in \mathbb{F}}\left(\pi\left(g\left(\begin{smallmatrix}
  c & t \\
  -1 & 0 \\
\end{smallmatrix}\right)\right)D_\infty\right)(1) \notag\\
&= \sum_{c \in \mathbb{F}}
(\pi(g)\pi(u_c)\pi(w)\pi(h_{t^{-1}})D_\infty)(1)
=\sum_{c \in \mathbb{F}} (\pi(g)\pi(u_c)\pi(w)\mu^{-1}(t^{-1})D_\infty)(1)\notag \\
&= \sum_{c \in \mathbb{F}} (\pi(g)\pi(u_c)(-q^{-1} \ve(\pi,
\mu^{-1},
\psi)\mu(t))(\mu+D_0))(1) \notag\\
&= -q^{-1} \ve(\pi, \mu^{-1}, \psi)\mu(t) \pi(g)\left[ \sum_{c\in
\mathbb{F}^\times} (q-1)^{-1} \sum_{\beta \in \hf\tx} \beta
\mu^{-1}(c) \gG(\mu\beta^{-1}, \psi)\beta + \mu + \sum_{c \in
\mathbb{F}} D_0
\right](1) \notag\\
&= -q^{-1} \ve(\pi, \mu^{-1}, \psi)\mu(t) \pi(g)\left[ (q-1)^{-1}
(q-1)\gG(1, \psi)\mu + \mu + qD_0\right](1) \notag\\
&= -\ve(\pi, \mu^{-1}, \psi) \mu(t) [\pi(g)D_0](1)= -\ve(\pi,
\mu^{-1}, \psi) \mu(t) \bW_{D_0}(g)\notag
\end{align}
for all $g \in G$. With respect to the basis $\{\bW_{D_0},
\bW_{D_\infty}\}$, the operator $T_{U_t}$ is represented by the
matrix
\[
{\footnotesize\left(
\begin{matrix}
  0 & ~-\ve(\pi, \mu^{-1}, \psi)\mu(t) \\
  -\ve(\pi, \mu, \psi)\mu(t^{-1}) & 0 \\
\end{matrix}
\right).}
\]
As
\[
\ve(\pi, \mu, \psi)\ve(\pi, \mu^{-1}, \psi) = \gG(\mu^2, \psi)
\gG(1, \psi) \gG(1, \psi) \gG(\mu^{-2}, \psi) = q \mu^2(-1) = q,
\]
the eigenvalues are $\pm\sqrt{q}$ with corresponding
eigenfunctions $\pm \sqrt{q}\ve(\pi, \mu, \psi)^{-1}\mu(t)
\bW_{D_0} - \bW_{D_\infty}$.

We have shown

\begin{thm}
For $t \in \mathbb{F}\tx$, the Cayley graph $X_{U_t} =
\mathrm{Cay}(G/U, U_t/U)$ is a $q$-regular Ramanujan graph with
the following eigenvalues$:$ $\mu(t)q$ of multiplicity one,
$-\mu(t)$ of multiplicity $q$ for $\mu^2 = 1$, and $\pm \sqrt q$
of multiplicity $(q+1)(q-3)/2$.
\end{thm}

Therefore the graph $X_{U_t}$ is bipartite for $t$ nonsquare since
$\mu(t) = -1$ for $\mu$ of order 2.  If $t$ is a square, the
graph has two connected components, one with square determinants
and one with nonsquares. In case $t = 1$, one component is
isomorphic to the graph $\Cay(\PSL_2(\mathbb{F})/U, UwU/U)$,
which, for the case $q = p$ a prime, is a cover of the Ramanujan
graph on cusps of the principal congruence subgroup $\Gamma(p)$ of
$\SL_2(\mathbb Z)$ studied by Gunnells in \cite{G04}. When $q =
p^e$ is a power of the prime $p$, a similar interpretation holds,
as we explain below.

The cusps of the Drinfeld upper half plane attached to the
rational function field $\mathbb F_p(T)$ and the (degree one)
place at infinity is represented by $\mathbb F_p(T) \cup
\{\infty\} = \mathbb{P}^1 (\mathbb F_p(T))$, on which the group
$\Gamma = \GL_2(\mathbb F_p[T])$ acts via fractional linear
transformations.  Endow a graph structure, called $X$, on
$\mathbb{P}^1 (\mathbb F_p(T))$ by defining two cusps $u, v$ to be
adjacent if there exists $\gamma \in \Gamma$ sending the cusp 0 to
$u$ and the cusp $\infty$ to $v$. Let $P$ be an irreducible
polynomial of degree $e$ over $\mathbb F_p$, and let $\Gamma(P)$
be the principal congruence subgroup of $\Gamma$ consisting of
matrices congruent to the identity matrix modulo $P$. Denote by
$X_P$ the quotient graph $\Gamma(P)\backslash X$. This is a
positive characteristic analog of the graph $G(p)$ in Gunnells's
paper \cite{G04}.

We proceed to show the connection between $X_P$ and
$\Cay(\PSL_2(\mathbb{F})/U, UwU/U)$. First observe the following
lemma, which can be deduced by the same argument as in the proof
of Lemma 1.42 in Shimura's book \cite{Sh94}.

\begin{lemma}\label{shi}
Let $t = \frac {a}{b}$ and $t' = \frac {c}{d}$ be two cusps in
$\mathbb{P}^1 (\mathbb F_p(T))$, where $a, b, c, d \in \mathbb
F_p[T]$ and $\gcd(a, b) = \gcd(c, d) = 1$. Then $t$ and $t'$ are
equivalent under $\Gamma(P)$ if and only if
$\left(\begin{smallmatrix}
  a \\
  b \\
\end{smallmatrix}\right) \equiv \lambda \left(\begin{smallmatrix}
  c \\
  d \\
\end{smallmatrix}\right) \mod P$
for some $\lambda \in \mathbb F_p^\times$.
\end{lemma}

Note that $\mathbb F_p[T]$ modulo $P$ is a finite field with
$p^e=q$ elements, which we identify with $\mathbb F$. By the lemma
above, the vertices in $X_P$ may be expressed as the column
vectors $\left(\begin{smallmatrix}
  a \\
  b \\
\end{smallmatrix}\right)$ modulo scalar multiplications by $\mathbb F_p^\times$,
where $a, b \in \mathbb F$, not both zero. The cusp $0$ is
represented by $\left(\begin{smallmatrix}
  0 \\
  1 \\
\end{smallmatrix}\right)$ and cusp $\infty$ by $\left(\begin{smallmatrix}
  1 \\
  0 \\
\end{smallmatrix}\right)$. The image of $\Gamma$ modulo $P$ is
\[
G' = \{ \gamma \in \GL_2(\mathbb F) : \det \gamma \in \mathbb
F_p^\times \}
\]
and the action of $\Gamma$ on the cusps of $\Gamma(P)$ becomes
(matrix) left multiplication by $G'$.

Observe that the vertices of $X_P$ are also the first (or second)
columns of the matrices in $G'$ modulo $\mathbb F_p^\times$. Since
\[
G' = \SL_2(\mathbb F)\cdot\{\left(\begin{smallmatrix}
  a & 0\\
  0 & 1 \\
\end{smallmatrix}\right) : a \in \mathbb F_p^\times \}
\]
acts transitively on the vertices of $X_P$ and the stabilizer of
$\left(\begin{smallmatrix}
  1 \\
  0 \\
\end{smallmatrix}\right)$ in $G'$ is the Borel subgroup of $G'$,
 the vertices of $X_P$ are represented by $\PSL_2(\mathbb F)/B$, where
$B$ denotes the Borel subgroup of $\PSL_2(\mathbb F)$. Suppose
that vertices $u, v$ of $X_P$ are adjacent, that is, there is a
matrix $\gamma \in G'$ such that $\gamma \left(\begin{smallmatrix}
  0 \\
  1 \\
\end{smallmatrix}\right) = u$ and $\gamma \left(\begin{smallmatrix}
  1 \\
  0 \\
\end{smallmatrix}\right) = v$. Then $v$ and $u$ are the first
and second columns of $\gamma$. In particular, this shows that the
neighbors of $\left(\begin{smallmatrix}
  1 \\
  0 \\
\end{smallmatrix}\right)$ are the second columns of the elements in the unipotent subgroup
$U$ of $G'$, which is also the unipotent subgroup of
$\PSL_2(\mathbb F)$. Note that they are the first columns of the
coset representatives of the double coset $UwU = \cup_{c \in
\mathbb F} \left(\begin{smallmatrix}
  c & ~-1\\
  1 & ~0 \\
\end{smallmatrix}\right)U.$ As the edge structure on $X_P$ is defined
by transporting the edges out of the vertex
$\left(\begin{smallmatrix}
  1 \\
  0 \\
\end{smallmatrix}\right)$ to other vertices using
left multiplications by $G'$, we have shown

\begin{prop} The graph $X_P$ on cusps of $\Gamma(P)$ is a quotient of  $\Cay(\PSL_2(\mathbb{F})/U, UwU/U)$, hence also a quotient of
$\Cay(G/U, UwU/U)$. Consequently, it is a $q$-regular Ramanujan
graph.
\end{prop}

\bigskip

\section{The $A$-graphs}

In this section, the group $H$ is $A$. We start with the space
$\cl_{\psi^0}(G/A)$. The group $G=UA \cup
UA\left(\begin{smallmatrix}
  0 & 1 \\
  -1 & 0 \\
\end{smallmatrix}\right)A \cup UA\left(\begin{smallmatrix}
  0 & 1 \\
  -1 & 1 \\
\end{smallmatrix}\right)A$ is a disjoint union of three double
cosets. For $\mu \ne 1$, one can easily see that the right
$A$-invariant space of $\Ind \mu$ is 1-dimensional, generated by
\[
f_\mu\left(\left(\begin{smallmatrix}
  y & x \\
  0 & 1 \\
\end{smallmatrix}\right)\left(\begin{smallmatrix}
  0 & 1 \\
  -1 & 1 \\
\end{smallmatrix}\right)A\right) = \mu(y) \; \text{for all $y \in
\mathbb{F}^\times, x \in \mathbb{F}$} \; \text{and} \; f_\mu(UA) =
f_\mu \left(UA\left(\begin{smallmatrix}
  0 & 1 \\
  -1 & 0 \\
\end{smallmatrix}\right)A\right) =0.
\]
If $\mu =1$, then the right $A$-invariant space is $3$-dimensional
generated by $\{f_1, f_2, f_3\}$ where $f_1, f_2$ and $f_3$ are
the characteristic functions of $UA, UA\left(\begin{smallmatrix}
  0 & 1 \\
  -1 & 1 \\
\end{smallmatrix}\right)A$ and $UA\left(\begin{smallmatrix}
  0 & 1 \\
  -1 & 0 \\
\end{smallmatrix}\right)A$, respectively. Clearly, $\cl_{\psi^0}(1 \circ \det, G/A)$ is
1-dimensional. So $\cl_{\psi^0}(\pi_1, G/A)$ is $2$-dimensional.
Next we fix a nontrivial additive character $\psi$ of
$\mathbb{F}$. For each irreducible representation $\pi$ of $G$ of
degree greater than 1, the function
\[
W_v(g) := (q-1)^{-1}\sum_{r \in \mathbb{F}^\times} (\pi(gh_r)v)(1)
\]
is right $A$-invariant for all $v \in \ck_\psi(\pi)$. From the
action of $\pi(h_r)$ on $\gh \in \hf^\times, D_0$ and $D_\infty$
described in \S 2, we find
\[
W_\gh(g) = \left\{\begin{array}{ll}
    (\pi(g)\gh)(1), & \hbox{if $\gh = 1$;} \\
    0, & \text{if $\gh \ne 1$,} \\
\end{array}\right. , \quad
W_{D_0}(g) = \left\{ \begin{array}{ll}
    (\pi(g)D_0)(1), & \hbox{if $\pi = \pi_1$;} \\
    0, & \hbox{if $\pi = \pi_\mu$ with $\mu \ne 1$,} \\
\end{array}\right.
\]
and $W_{D_\infty}(g) = 0$ for all $g \in G$. Since $\dim \cl(G/A)
= q(q-1)$, we conclude by dimension counting that for a principal
series representation $\pi_\mu$, the Steinberg representation
$\pi_\mu$ with $\mu$ quadratic, and for a discrete series
representation $\pi_\gL$, the space of right $A$-invariant vectors
in its Whittaker model is 1-dimensional generated by $\{W_1\}$,
and the space of right $A$-invariant vectors is 2-dimensional
generated by $\{W_1, W_{D_0}\}$ for the Steinberg representation
$\pi_1$.

We have proven

\begin{prop}
\begin{itemize}
    \item[(1)] The space $\cl_{\psi^0}(\pi_\mu, G/A)$ is $1$-dimensional
if $\mu \ne 1$ and $2$-dimensional if $\mu = 1$;
 and the space $\cl_{\psi^0}(1 \circ \det, G/A)$ is $1$-dimensional.
    \item[(2)] For nontrivial $\psi$'s, the space $\cl_\psi(\pi, G/A)$ is $2$-dimensional if $\pi$ is the Steinberg
               representation $\pi_1$, and $1$-dimensional otherwise.
\end{itemize}
\end{prop}

The group $G$ has $q+4$ $A$-double cosets. Except for
$A\left(\begin{smallmatrix}
  1 & 0 \\
  0 & 1 \\
\end{smallmatrix}\right)A = A$, $A\left(\begin{smallmatrix}
  0 & 1 \\
  -1 & 0 \\
\end{smallmatrix}\right)A = \left(\begin{smallmatrix}
  0 & 1 \\
  -1 & 0 \\
\end{smallmatrix}\right)A$, each of the the remaining $q+2$ double cosets $A\left(\begin{smallmatrix}
  0 & 1 \\
  -1 & 1 \\
\end{smallmatrix}\right)A$, $A\left(\begin{smallmatrix}
  1 & 1 \\
  0 & 1 \\
\end{smallmatrix}\right)A$ and $A\left(\begin{smallmatrix}
  1 & \gd-c \\
  1 & 1-c \\
\end{smallmatrix}\right)A$, $c \in \mathbb{F}$, is the disjoint union of
$q$ A-cosets. Write $A_c = A\left(\begin{smallmatrix}
  1 & \gd-c \\
  1 & 1-c \\
\end{smallmatrix}\right)A$ and $A_\nf = A\left(\begin{smallmatrix}
  1 & 1 \\
  0 & 1 \\
\end{smallmatrix}\right)A$ for short. Among these, we rule out $A_\delta$
and $A_\infty$ since they are contained in conjugates of the Borel
subgroup of $G$, and $A_{1}$ and $A\left(\begin{smallmatrix}
  0 & 1 \\
  -1 & 1 \\
\end{smallmatrix}\right)A$ since they are not their own inverse. The
remaining symmetric double cosets $A_c, c\in \mathbb{F}
\smallsetminus \{1, \delta\}$ are of interest to us. For such $c$,
we express $A_c$ as the union of $A$-cosets:
\[
A_c = \bigcup_{x \in \mathbb{F}\tx} \left(\begin{smallmatrix}
  x & x(\gd-c) \\
  1 & 1-c \\
\end{smallmatrix}\right)A.
\]
Our computations on the space $\cl_{\psi^0}(G/A)$ will use the
following facts repeatedly:
\begin{itemize}
\item[(i)] $\left(\begin{smallmatrix}
  x & x(\gd-c) \\
  1 & 1-c \\
\end{smallmatrix}\right)A = \left(\begin{smallmatrix}
  \frac{x(\gd-1)}{1-c} & x \\
  0 & 1 \\
\end{smallmatrix}\right)\left(\begin{smallmatrix}
  0 & 1 \\
  -1 & 1 \\
\end{smallmatrix}\right)A$ and $\left(\begin{smallmatrix}
  0 & 1 \\
  -1 & 0 \\
\end{smallmatrix}\right)\left(\begin{smallmatrix}
  x & x(\gd-c) \\
  1 & 1-c \\
\end{smallmatrix}\right)A = \left(\begin{smallmatrix}
  \frac{1-\gd}{x(\gd-c)} & -\frac{\gd}{x} \\
  0 & 1 \\
\end{smallmatrix}\right)\left(\begin{smallmatrix}
  0 & 1 \\
  -1 & 1 \\
\end{smallmatrix}\right)A$ for all $x \in \mathbb{F}^\times$,
\item[(ii)] $\left(\begin{smallmatrix}
  0 & 1 \\
  -1 & 1 \\
\end{smallmatrix}\right)\left(\begin{smallmatrix}
  x & x(\gd -c) \\
  1 & 1-c \\
\end{smallmatrix}\right)A = \left(\begin{smallmatrix}
  \frac{x(\gd-1)}{(x-1)(x(\gd-c)-(1-c))} & -\frac{1}{1-x} \\
  0 & 1 \\
\end{smallmatrix}\right)\left(\begin{smallmatrix}
  0 & 1 \\
  -1 & 1 \\
\end{smallmatrix}\right)A$ for all $x \in \mathbb{F}^\times, x \ne 1,
\frac{c-1}{c-\gd}$,  \\ $\left(\begin{smallmatrix}
  0 & 1\\
  -1 & 1 \\
\end{smallmatrix}\right)\left(\begin{smallmatrix}
  1 & \gd-c \\
  1 & 1-c \\
\end{smallmatrix}\right)A = \left(\begin{smallmatrix}
  \frac{1}{\gd-1} & \frac{c-1}{\gd-1} \\
  0 & 1 \\
\end{smallmatrix}\right)A$ and $\left(\begin{smallmatrix}
  0 & 1 \\
  -1 & 1 \\
\end{smallmatrix}\right)\left(\begin{smallmatrix}
  \frac{1-c}{\gd-c} & 1-c \\
  1 & 1-c \\
\end{smallmatrix}\right)A = \left(\begin{smallmatrix}
  \frac{(1-c)(\gd-c)}{\gd-1} & \frac{\gd-c}{\gd-1} \\
  0 & 1 \\
\end{smallmatrix}\right)\left(\begin{smallmatrix}
  0 & 1 \\
  -1 & 0 \\
\end{smallmatrix}\right)A$.
\end{itemize}
For $\mu \ne 1$, $\cl_{\psi^0}(\pi_\mu, G/A)$ is 1-dimensional
generated by $\{f_\mu\}$ with
$f_\mu\left(\left(\begin{smallmatrix}
  0 & 1 \\
  -1 & 1 \\
\end{smallmatrix}\right)\right) = 1$, so it is an eigenspace of $T_{A_c}$ with the
eigenvalue
\[
{\small\begin{aligned} &\gl =
(T_{A_c}f_\mu)\left(\left(\begin{smallmatrix}
  0 & 1 \\
  -1 & 1 \\
\end{smallmatrix}\right)\right) = \sum_{x \in \mathbb{F}^\times} f_\mu\left(\left(\begin{smallmatrix}
  0 & 1 \\
  -1 & 1 \\
\end{smallmatrix}\right)\left(\begin{smallmatrix}
  x & x(\gd-c) \\
  1 & 1-c \\
\end{smallmatrix}\right)\right) \\
&= f_\mu\left(\left(\begin{smallmatrix}
  \frac{1}{\gd-1} & \frac{c-1}{\gd-1} \\
  0 & 1 \\
\end{smallmatrix}\right)\right) + f_\mu\left(\left(\begin{smallmatrix}
  \frac{(1-c)(\gd-c)}{\gd-1} & \frac{\gd-c}{\gd-1} \\
  0 & 1 \\
\end{smallmatrix}\right)\left(\begin{smallmatrix}
  0 & 1 \\
  -1 & 0 \\
\end{smallmatrix}\right)\right) + \sum_{x \in \mathbb{F}^\times \atop x \ne 1, \frac{1-c}{\gd-c}}
f_\mu\left(\left(\begin{smallmatrix}
  \frac{x(\gd-1)}{(x-1)(x(\gd-c)-(1-c))} & -\frac{1}{x-1} \\
  0 & 1 \\
\end{smallmatrix}\right)\left(\begin{smallmatrix}
  0 & 1 \\
  -1 & 1 \\
\end{smallmatrix}\right)\right) \\
&= \sum_{x \in \mathbb{F}^\times \atop x \ne 1, \frac{1-c}{\gd-c}}
\mu\left( {\small\tfrac{x(\gd-1)}{(x-1)(x(\gd-c)-(1-c))}}\right).
\end{aligned}}
\]
Therefore $|\gl| \le 2\sqrt{q}$ by Theorem 3 in Chapter 6 of
\cite{L96}. If $\mu = 1$, then $\cl_{\psi^0}(\pi_1, G/A) \oplus
\cl_{\psi^0}(1 \circ \det, G/A)$ is 3-dimensional, generated by
$\{f_1, f_2, f_3\}$. For $i \in \{1, 2, 3 \}$, we have
\[
\begin{split}
&(T_{A_c} f_i)(UA) = \sum_{x \in \mathbb{F}^\times}
f_i\left(\left(\begin{smallmatrix}
  x & x(\gd-c)\\
  1 & 1-c\\
\end{smallmatrix}\right)\right) = \sum_{x \in \mathbb{F}^\times} f_i\left(\left(\begin{smallmatrix}
  \frac{x(\gd-1)}{1-c} & x\\
  0 & 1 \\
\end{smallmatrix}\right)\left(\begin{smallmatrix}
  0 & 1\\
  -1 & 1 \\
\end{smallmatrix}\right)\right), \\
&(T_{A_c} f_i)(UA\left(\begin{smallmatrix}
  0 & 1\\
  -1 & 1 \\
\end{smallmatrix}\right)A) = \sum_{x \in \mathbb{F}^\times}
f_i\left(\left(\begin{smallmatrix}
  0 & 1\\
  -1 & 1 \\
\end{smallmatrix}\right)\left(\begin{smallmatrix}
  x & x(\gd-c)\\
  1 & 1-\gd \\
\end{smallmatrix}\right)\right) \\
& \,\, =f_i\left(\left(\begin{smallmatrix}
  \frac{1}{1-\gd} & \frac{c-1}{\gd-1} \\
  0 & 1 \\
\end{smallmatrix}\right)\right) + f_i\left(\left(\begin{smallmatrix}
  \frac{(1-c)(\gd-c)}{\gd-1} & \frac{\gd-c}{\gd-1} \\
  0 & 1 \\
\end{smallmatrix}\right)\left(\begin{smallmatrix}
  0 & 1 \\
  -1 & 0 \\
\end{smallmatrix}\right)\right) +\sum_{x \in \mathbb{F}^\times \atop a \ne 1,
\frac{1-c}{\gd-c}} f_i\left(\left(\begin{smallmatrix}
 \frac{x(\gd-1)}{(x-1)(x(\gd-c)-(1-c))} & -\frac{\gd}{x-1}\\
  0 & 1 \\
\end{smallmatrix}\right)\left(\begin{smallmatrix}
  0 & 1\\
  -1 & 1 \\
\end{smallmatrix}\right)\right), \\
&(T_{A_c} f_i)(UA\left(\begin{smallmatrix}
  0 & \gd\\
  1 & 0 \\
\end{smallmatrix}\right)A) = \sum_{x \in \mathbb{F}^\times}
f_i\left(\left(\begin{smallmatrix}
  0 & 1\\
  -1 & 0 \\
\end{smallmatrix}\right)\left(\begin{smallmatrix}
  x & x(\gd-c)\\
  1 & 1-c \\
\end{smallmatrix}\right)\right) = \sum_{x \in \mathbb{F}^\times}f_i\left(\left(\begin{smallmatrix}
  \frac{\gd-1}{x(\gd-c)} & -\frac{1}{x} \\
  0 & 1 \\
\end{smallmatrix}\right)\left(\begin{smallmatrix}
  0 & 1\\
  -1 & 1 \\
\end{smallmatrix}\right)\right).
\end{split}
\]
This shows that
\[
(T_{A_c} f_1) = f_2, (T_{A_c} f_2) = (q-1)f_1+(q-3)f_2+(q-1)f_3
\quad \text{and} \quad (T_{A_c} f_3) = f_2.
\]
So with respect to the basis $f_1, f_2, f_3$, the operator
$T_{A_c}$ is represented by the matrix
\[
{\footnotesize \left(\begin{matrix}
  0 & q-1 & 0\\
  1 & q-3 & 1\\
  0 & q-1 & 0
\end{matrix}\right).}
\]
Clearly, $0$ is an eigenvalue of $T_{A_c}$. Further,
$\cl_{\psi^0}(1 \circ \det, G/A)$ is a 1-dimensional eigenspace
with the eigenvalue $q-1$. Hence $0, q-1$ and $-2$ (from trace
computation) are eigenvalues of $T_{A_c}$ with the eigenfunctions
$f_1 - f_3, f_1+f_2+f_3$ and $\frac{1-q}{2}f_1 + f_2 +
\frac{1-q}{2}f_3$, respectively.

Now we turn to a nontrivial additive character $\psi$. Recall that
$\cl_\psi (\pi, G/A)$ is generated by $\{W_1\}$ for $\pi \ne
\pi_1$ and by $\{W_1, W_{D_0} \}$ for $\pi = \pi_1$. We shall
compute the actions of $T_{A_c}$ on these two functions. Note that
\[
\left(\begin{smallmatrix}
  x & x(\gd-c) \\
  1 & 1-c \\
\end{smallmatrix}\right) = \left(\begin{smallmatrix}
  x & 0 \\
  0 & 1 \\
\end{smallmatrix}\right)\left(\begin{smallmatrix}
  1 & 1 \\
  0 & 1 \\
\end{smallmatrix}\right)\left(\begin{smallmatrix}
  0 & 1 \\
  -1 & 0 \\
\end{smallmatrix}\right)\left(\begin{smallmatrix}
  1 & (1-c)(1-\gd)^{-1} \\
  0 & 1 \\
\end{smallmatrix}\right)\left(\begin{smallmatrix}
  (1-\gd)^{-1} & 0 \\
  0 & 1\\
\end{smallmatrix}\right)
\left(\begin{smallmatrix}
  \gd-1 & 0 \\
  0 & \gd-1 \\
\end{smallmatrix}\right) \quad \text{for all $x \in \mathbb{F}^\times$}.
\]
Using the Kirillov model of $\pi$ given in \S 2, we compute first
\[
{\small\begin{aligned} &(T_{A_c} W_1)(g) = \sum_{x \in
\mathbb{F}^\times} W_1\left(g \left(\begin{smallmatrix}
  x & x(\gd-c) \\
  1 & 1-c \\
\end{smallmatrix}\right)\right)
= \sum_{x} \left[\pi\left(g \left(\begin{smallmatrix}
  x & x(\gd-c) \\
  1 & 1-c \\
\end{smallmatrix}\right)\right)1\right](1) \\
&= \sum_x [\pi(g) \pi(h_x) \pi(u_1) \pi(w)
\pi(u_{(1-c)(1-\gd)^{-1}})
\pi(h_{(1-\gd)^{-1}})1](1) \\
&= (q-1)^{-2}q^{-1} \sum_{\beta \in \hf^\times}
\beta((1-c)(1-\gd)^{-1}) \gG(\beta^{-1}, \psi)\ve(\pi, \beta,
\psi)\sum_{\gamma \in \hf^\times}
\gG(\beta^{-1}\gamma^{-1})[\pi(g)\sum_{x} \gamma(x)
\gamma](1) \\
& \quad +(q-1)^{-1} \sum_{\beta \in \hf^\times}
\beta((1-c)(1-\delta)^{-1})\gG(\beta^{-1}, \psi) [\pi(g)\sum_{x} \pi(h_x)e(\beta)](1) \\
&= (q-1)^{-1}q^{-1} \sum_{\beta} \beta((1-c)(1-\gd)^{-1})
\gG(\beta^{-1}, \psi)^2\ve(\pi, \beta, \psi)[\pi(g)1](1) \\
& \quad +(q-1)^{-1} \sum_{\beta} \beta((1-c)(1-\delta)^{-1})
\gG(\beta^{-1}, \psi) [\pi(g)\sum_{x} \pi(h_x)e(\beta)](1).
\end{aligned}}
\]
The sum $\sum_{x \in \mathbb{F}^\times} \pi(h_x) e(\beta) =
-q^{-1}(q-1)\ve(\pi, \beta, \psi)(q^2 -1)\gd_{\beta,1}D_0$ if
$\pi$ is the Steinberg representation $\pi_1$, and 0 otherwise. In
case $\pi$ is the Steinberg representation $\pi_1$, we also need
to find the action of $T_{A_c}$ on $W_{D_0}$, which is, for $g \in
G$,
\[{\small\begin{aligned} (T_{A_c} W_{D_0})(g) &= \sum_{x \in
\mathbb{F}^\times} W_{D_0}\left(g \left(\begin{smallmatrix}
  x & x(\gd-c) \\
  1 & 1-c \\
\end{smallmatrix}\right)\right)
= \sum_{x} \left[\pi\left(g \left(\begin{smallmatrix}
  x & x(\gd-c) \\
  0 & 1-c \\
\end{smallmatrix}\right)\right)D_0\right](1) \\
&= \sum_x [\pi(g) \pi(h_x) \pi(u_1) \pi(w)
\pi(u_{(1-c)(1-\gd)^{-1}})\pi(h_{(1-\gd)^{-1}}) D_0](1) \\
&= \sum_x [\pi(g)\pi(h_{x^{-1}})\pi(u_1)\pi(w)
\pi(u_{(1-c)(1-\gd)^{-1}})D_0](1) \\
&= \sum_x \pi(g)\pi(h_{x^{-1}})\pi(u_1) [-q^{-1}
\ve(\pi_1, 1, \psi)(1+D_0)](1) \\
&= -q^{-1} \sum_x \pi(g)\pi(h_{x^{-1}})
\left[(q-1)^{-1}\sum_{\beta \in \hf\tx} \gG(\beta^{-1}, \psi)\beta + D_0\right](1) \\
&= -q^{-1}\pi(g)\left[ (q-1)^{-1} \sum_\beta \gG(\beta^{-1}, \psi)
\sum_x \beta(x) \beta + \sum_x  D_0 \right](1) \\
&= q^{-1}[\pi(g)1](1) - q^{-1}(q-1)[\pi(g)D_0](1) = q^{-1}W_1(g) -
q^{-1}(q-1)W_{D_0}(g)
\end{aligned}}
\]
since $\ve(\pi_1, 1, \psi) = 1$.

To compute the eigenvalues and the eigenfunctions of $T_{A_c}$, we
begin with the 2-dimensional space $\cl_\psi(\pi_1, G/A)$
generated by $\{W_1, W_{D_0}\}$. The identities
\[
\ve(\pi_1, \beta, \psi) = \gG(\beta, \psi)^2 \quad \text{and}
\quad \gG(\beta,\psi) \gG(\beta^{-1},\psi) = \left\{
\begin{array}{ll}
    1, & \hbox{if $\beta =1$;} \\
    q, & \hbox{if $\beta \ne 1$,} \\
\end{array}
\right.
\]
yield the expression
\[
(T_{A_c}W_1)(g) = -(q-1)^{-1}q^{-1}(q^2-1)W_1(g) +
q^{-1}(q^2-1)W_{D_0}(g).
\]
Hence with respect to the basis $\{W_1, W_{D_0}\}$, the operator
$T_{A_c}$ is represented by the matrix
\[
{\footnotesize \left(\begin{matrix}
  -\frac{q+1}{q} & ~\frac{1}{q} \\
  \frac{q^2-1}{q} & ~-\frac{q-1}{q} \\
\end{matrix}\right).}
\]
Thus $0$ and $-2$ are the eigenvalues of $T_{A_c}$ with
eigenfunctions $W_1+(q+1)W_{D_0}$ and $(q-1)W_1 - W_{D_0}$,
respectively.

Next we consider the case where $\pi \ne \pi_1$. Then
$\cl_\psi(\pi, G/A)$ is 1-dimensional, so it is an eigenspace
generated by $\{W_1\}$ with the eigenvalue
\[
\gl = (q-1)^{-1}q^{-1}\sum_{\beta \in \hf^\times}
\beta((1-c)(1-\gd)^{-1}) \gG(\beta^{-1}, \psi)^2\ve(\pi, \beta,
\psi).
\]
When $\pi$ is a discrete series representation $\pi_\gL$, we have
$\ve(\pi, \beta, \psi) =  - \sum_{z \in
\mathbb{\mathbb{E}}^\times} \Lambda(z) \beta(\nr \, z)\psi(\tr \,
z)$ so that
\[
{\small\begin{aligned} &\sum_{\beta \in \hf^\times}
\beta((1-c)(1-\gd)^{-1})
\gG(\beta^{-1},\psi)^2\ve(\pi, \beta, \psi) \\
&= - \sum_{\beta \in \hf^\times} \beta((1-c)(1-\gd)^{-1})\sum_{x
\in \mathbb{F}^\times} \beta^{-1}(x)\psi(x) \sum_{y \in
\mathbb{F}^\times} \beta^{-1}(y) \psi(y) \sum_{z \in
\mathbb{\mathbb{E}}^\times} \gL(z) \beta(\nr
\, z) \psi(\tr \, z) \\
&= -(q-1) \sum_{x \in \mathbb{F}^\times} \sum_{z \in
\mathbb{\mathbb{E}}^\times} \psi(x)
\psi\left(\tfrac{(1-c)\nr\, z}{(1-\gd)x}\right)\gL(z) \psi(\tr\,z) \\
&= -(q-1) \sum_{z \in \mathbb{\mathbb{E}}^\times} \gL(z)
\psi(\tr\,z) \sum_{x \in \mathbb{F}^\times} \psi\left(x
+\tfrac{(1-c)\nr\, z}{(1-\gd)x}\right).
\end{aligned}}
\]
Thus $|\gl| \le (q-1)^{-1}q^{-1}(q-1)q(2\sqrt{q}) = 2 \sqrt{q}$ by
Corollary 5 in Chapter 6 of \cite{L96}.

When $\pi$ is a principal series representation or a Steinberg
representation $\pi_\mu$ with $\mu \ne 1$, we have $\ve(\pi,
\beta, \psi) = \gG(\mu\beta, \psi) \gG(\mu^{-1}\beta, \psi)$ so
that
\[
{\small\begin{aligned} &\sum_{\beta \in \hf^\times}
\beta((1-c)(1-\gd)^{-1}) \gG(\beta^{-1},
\psi)^2\ve(\pi, \beta, \psi) \\
&=  \sum_{\beta \in \hf^\times} \beta((1-c)(1-\gd)^{-1})\sum_{x
\in \mathbb{F}^\times} \beta^{-1}(x)\psi(x) \sum_{y \in
\mathbb{F}^\times} \beta^{-1}(y) \psi(y) \sum_{u \in
\mathbb{F}^\times} \mu\beta(u)\psi(u)
\sum_{z \in \mathbb{F}^\times} \mu^{-1}\beta(z) \psi(z)\\
&= (q-1) \sum_{x, u, z \in \mathbb{F}^\times} \psi(x) \psi \left(
\tfrac{uz(1-c)}{(1-\gd)x} \right) \mu(u) \psi(u) \mu^{-1}(z) \psi(z) \\
&= (q-1) \sum_{u \in \mathbb{F}^\times} \mu(u)\psi(u) \sum_{z \in
\mathbb{F}^\times} \mu^{-1}(z) \psi(z) \sum_{x \in
\mathbb{F}^\times} \psi\left( x + \tfrac{uz(1-c)}{(1-\gd)x}
\right).
\end{aligned}}
\]
Thus $|\gl| \le (q-1)^{-1}q^{-1}(q-1)\sqrt{q} \sqrt{q} (2\sqrt{q})
= 2 \sqrt{q}$ by the same corollary.

\medskip
We have found all eigenvalues and eigenfunctions of the operator
$T_{A_c}$ on the space $\cl(G/A)$. The eigenvalues are the
spectrum of the Cayley graph $\Cay(G/A, A_c/A)$; the estimates
indicate that it is almost a Ramanujan graph. We record this in

\begin{thm}
For $c \ne 1, \gd$, all nontrivial eigenvalues of the Cayley graph
$X_{A_c} = \mathrm{Cay}(G/A, A_c/A)$ have absolute value at most
$2 \sqrt q$.
\end{thm}


\begin{thebibliography}{99}

\bibitem{AC92} J. Angel, N. Celniker, S. Poulos, A.
Terras, C. Trimble and E. Velasquez, Special functions on finite
upper half planes, Contemp. Math. 138 (1992) 1-26.

\bibitem{CL03} C.-L. Chai and W.-C. W. Li,  Character sums, automorphic forms,
equidistribution, and Ramanujan graphs, Part I. The Kloosterman
sum conjecture over function fields, Forum Math. 15 Issue 5 (2003)
679-699.

\bibitem{CL04} \bysame, Character sums, automorphic forms,
equidistribution, and Ramanujan graphs, Part II. Eigenvalues of
Terras graphs, Forum Math. 16 (2004) 631-661.

\bibitem{CP93} N. Celniker, S. Poulos, A.
Terras, C. Trimble and E. Velasquez, Is there life on finite upper
half planes?, Contemp. Math. 143 (1993) 65-88.

\bibitem{De04} M. Dedeo, D. Lanphier and M. Minei, The spectrum of Platonic
graphs over finite fields, preprint, 2004.

\bibitem{Dr88} V. G. Drinfeld, The proof of Petersson's conjecture
for $\GL(2)$ over a global field of characteristic $p$, Functional
Anal. Appl. 22 (1988) 28-43.

\bibitem{Ei57} M. Eichler, Eine Verallgemeinerung der Abelschen Integrale, Math. Z. 67 (1957) 267--298.

\bibitem{G04} P. E. Gunnells, Some elementary Ramanujan graphs, 
Geometriae Dedicata, to appear.

\bibitem{L92} W.-C. W. Li, Character sums and abelian Ramanujan graphs, J.
Number Theory 41 (1992) 199-214.

\bibitem{L96} \bysame, Number Theory with
Applications, World Scientific, Singapore, 1996.

\bibitem{L99} \bysame, Eigenvalues of Ramanujan graphs, Emerging applications of number theory
(Minneapolis, MN, 1996),  387--403, IMA Vol. Math. Appl., 109,
Springer, New York, 1999.

\bibitem{L05} \bysame, Character sums over norm groups, Finite Fields Appl.,
to appear.

\bibitem{LSA83} W.-C. W. Li and J. Soto-Andrade, Barnes' identities and representations of
${\rm GL}(2)$. I. Finite field case., J. Reine Angew. Math. 344
(1983) 171-179.

\bibitem{LPS88} A. Lubotzky, R. Phillips and P. Sarnak, Ramanujan
graphs, Combinatorica 8 (1988) 261-277.

\bibitem{Ma88} G. Margulis, Explicit group theoretic constructions
of combinatorial schemes and their application to the design of
expanders and concentrators, J. Prob. of Info. Trans. (1988)
39-46.

\bibitem{Mo94} M. Morgenstern, Existence and explicit
constructions of $q+1$ regular Ramanujan graphs for every prime
power $q$, J. Comb. Theory, series B 62 (1994) 44-62.

\bibitem{PS83} I. I. Piatetski-Shapiro, Complex Representations of ${\rm GL}(2,K)$ for Finite Field
$K$, Contemporary Math. 16, Amer. Math. Soc., Providence, 1983.

\bibitem{Sh59} G. Shimura, Sur les int\'egrales attach\'ees aux formes
automorphes, J. Math. Soc. Japan 11 (1959) 291-311.

\bibitem{Sh94} \bysame, Introduction to the Arithmetic Theory of Automorphic Forms,
Publications of the Mathematical Society of Japan, vol. 11,
Princeton University Press, Princeton, NJ, 1994.

\bibitem{Ta99} A. Terras, Fourier Analysis on Finite Groups and Applications,
London Math Soc., Student Texts 43, Cambridge Univ. Press, Cambridge, UK, 1999.
\end{thebibliography}
\end{document}